\documentclass[11pt]{article}
\usepackage{amsmath,amsfonts,amssymb,mathrsfs,cases,mathdots,colordvi,url,extarrows,graphicx}
\usepackage{relsize}
\usepackage{mathrsfs}
\usepackage{color}
\usepackage{indentfirst}
\usepackage{float}

\newtheorem{defn}{Definition}[section]

\newtheorem{lemma}{Lemma}[section]
\newtheorem{thm}{Theorem}[section]
\newtheorem{prop}{Proposition}[section]
\newtheorem{cor}{Corollary}[section]
\newtheorem{example}{Example}[section]
\newtheorem{remark}{Remark}[section]

\renewcommand{\Box}{\rule{2.2mm}{2.2mm}}

\def\beginproof{\par\noindent {\bf Proof.}\ \ }
\def\endproof{\hskip .5cm $\Box$ \vskip .5cm}

\def\beginproof{\par\noindent {\bf Proof.}\ \ }
\def\endproof{\hskip .5cm $\Box$ \vskip .5cm}

\begin{document}
\title{ Directional derivative of the value function for  parametric set-constrained optimization problems}
\author{Kuang Bai\thanks{Department of Applied Mathematics, The Hong Kong Polytechnic University, Hong Kong, People's Republic of China. The research of this author was partially supported by the NSFC Grant 12201531 and the Hong Kong Research Grants Council PolyU153036/22p. Email: kuang.bai@polyu.edu.hk. } and Jane J. Ye\thanks{Corresponding author. Department of Mathematics and Statistics, University of Victoria, Canada. The research of this author was partially
supported by NSERC. Email: janeye@uvic.ca.}\\
{\small Dedicated to Boris Mordukhovich on the occasion of his seventy-fifth birthday}
}

\maketitle
\begin{abstract}

	This paper is concerned with the directional derivative of the value function for a very general set-constrained optimization problem under perturbation.
	Under reasonable assumptions, we obtain upper and lower estimates for the upper and lower Dini directional derivative of the value function respectively, from which we obtain Hadamard directional differentiability of the value function  when the set of multipliers is a singleton. Our results do not require convexity of the set involved. Even in the case of a parametric nonlinear program, 
	 our results improve the classical ones in that our regularity conditions are weaker and the directional solution set is used which is in general smaller than its nondirectional counterparts.
	
	\vskip 10 true pt
	

	\vskip 10 true pt
	

\end{abstract}

\newpage
\section{Introduction}
In this paper we consider the parametric set-constrained optimization problem in the form:
\begin{equation*} 
(P_x)~~~~\quad \min_{y\in\mathbb R^m} f(x,y)
\quad s.t.\  P(x,y)\in C,
\end{equation*}
where $x\in\mathbb R^n$ is the parameter or perturbation, $f:\mathbb R^{n+m}\rightarrow\mathbb R,\ P:\mathbb R^{n+m}\rightarrow\mathbb R^p$ are continuously differentiable, and $C\subseteq\mathbb R^p$ is a closed set. 
We define the feasible solution map as 
$${\cal F}(x):=\{y \in \mathbb R^m | P(x,y)\in C\},$$
  the (optimal) value function/marginal function as
$$V(x):=\inf_y \{f(x,y)| y\in {\cal F}(x)\}$$ 
and the optimal solution set as
$$S(x):=\arg\min_y \{f(x,y)| y\in {\cal F}(x)\}.$$

Recently we have studied the directional subdifferential of the value function in \cite{BYCOT}. The main goal of this paper is to study the directional derivative of the value function. 
In the case of fixed feasible set, i.e., $P(x,y)=y$, as it is commented by Bonnans and Shapiro in \cite[page 272]{BS},  the directionally differentiable behavior is typical for the value function. Indeed, by Danskin's theorem \cite{Danskin},  if $P(x,y)=y$ and $C$ is compact then $V(x)$ is directionally differentiable in any direction and the directional derivative of $V(x)$ at $\bar x$ in direction $u$ can be calculated by
$$V'(\bar x; u)=\min_{y\in S(\bar x)} \nabla_x f(\bar x, y)^T u.$$ Moreover by \cite[Theorem 4.13]{BS}, the compactness of set $C$ can be replaced by a weaker condition such as the inf-compactness condition.
As commented by Bonnans and Shapiro in \cite[page 278]{BS}, {\it ``It is considerably more difficult to investigate differentiability properties of the optimal value function in cases where the corresponding feasible set is also subject to perturbations".} In the case where $C=\{0\}_{p_1}\times  \mathbb{R}^{p_2}$ with $p_1+p_2=p$, the parametric optimization problem is a parametric nonlinear program. Under the uniform compactness {of the solution map $S( x)$ near $\bar x$} and the Mangasarian-Fromovitz constraint qualification (MFCQ)  at each $y\in S(\bar x)$, Gauvin and  Dubeau \cite[Corollary 4.3]{GD} obtained bounds of the upper/lower Dini directional derivative of the value function:
\begin{equation}\label{minmax}
\inf_{y\in S(\bar x)} \min_{\lambda \in \Lambda(\bar x, y)}\nabla_x L(\bar x, y,\lambda)^T u \leq V_-'(\bar x;u)\leq V_+'(\bar x;u)\leq \min_{y\in S(\bar x)} \max_{\lambda \in \Lambda(\bar x, y)}\nabla_x L(\bar x, y,\lambda)^T u,
\end{equation}{where $L(\bar x,y,\lambda):=f(\bar x,y)+P(\bar x,y)^T\lambda$ is the Lagrange function and $\Lambda(\bar x, y)$ is the set of Lagrange multipliers for problem $(P_{\bar x})$.
Moreover if the multiplier set $\Lambda(\bar x, y)=\{\lambda(y)\}$ is a singleton, then   the upper bound is equal to the lower bound and the value function is then directionally differentiable and 
$$V'(\bar x; u)=\inf_{y\in S(\bar x)}\nabla_x L(\bar x, y,\lambda(y))^T u. $$ These results are extended to cover the case where $C$ is an arbitrary closed convex set by Bonnans and Shapiro \cite[Theorem 4.26]{BS}.

In general there is a gap between the upper and lower bounds in (\ref{minmax}).
Under some additional assumptions the gap can be closed and one can conclude the directional differentiability and the formula
$$V'(\bar x; u)=\min_{y\in S(\bar x)}\max_{\lambda \in \Lambda(\bar x, y)}\nabla_x L(\bar x, y,\lambda)^T u. $$ For parametric nonlinear programs, 
 Janin \cite[Corollary 3.4]{Janin}  obtained such a result under both the constant rank constraint qualification 
 and MFCQ. For the parametric program in the form of $(P_x)$ with a closed convex set $C$, these types of  results are given by Bonnans and Shapiro in \cite[Theorems 4.24 and 4.25]{BS}.

In this paper, we aim at studying the directional derivative of the value function of the very general parametric program $(P_x)$. First we extend the result of Bonnans and Shapiro \cite[Theorem 4.26]{BS} to allow $C$ to be nonconvex. Moreover we obtain
the following  estimates  \begin{equation*}
\inf_{y\in S(\bar x;u)} \min_{\lambda \in \Lambda^c(\bar x, y)}\nabla_x L(\bar x, y,\lambda) u \leq V_-'(\bar x;u)\leq V_+'(\bar x;u)\leq \inf_{y\in S(\bar x;u)} \max_{\lambda \in \Lambda^c(\bar x, y)}\nabla_x L(\bar x, y,\lambda) u,
\end{equation*}
where $\Lambda^c(\bar x, y)$ denotes the set of Clarke multipliers for problem $(P_{\bar x})$ and 
$S(\bar x;u)$ is a subset of the solution set $S(\bar x)$ associated with the direction $u$. Furthermore if the set of Clarke multipliers  $\Lambda^c(\bar x, y)=\{\lambda(y) \}$ is a singleton for all $y\in S(\bar x;u)$, then the value function is directionally differentiable and 
$$V'(\bar x; u)=\inf_{y\in S(\bar x;u)}\nabla_x L(\bar x, y,\lambda(y)) u. $$
Our result has even improved the classical one 
\cite[Theorem 4.26]{BS} in the case when the set $C$ is convex since the regularity conditions we require are weaker than Robinson's constraint qualification and the set of directional solutions $S(\bar x;u)$ may be strictly included in  the solution set $S(\bar x)$. In particular our result for the parametric nonlinear program has improved the classical result of Gauvin and  Dubeau \cite[Corollary 4.3]{GD} in that our regularity condition is weaker than MFCQ.

{ The rest of this paper is organized as follows. In Section 2, we give some preliminaries and preliminary results used in the paper. In Section 3, we propose first-order sufficient conditions for directional Robinson stability. In Section 4, we develop upper-/lower-bound for Dini directional derivatives of the value function for parametric set-constrained problems. Finally, an example is used to illustrate the sharpness of our estimates.}

\section{Preliminaries and preliminary results} We first give notations that will be used in the paper.  Let $\Omega$ be a set. By $x^k\xrightarrow{\Omega}\bar{x}$ we mean $x^k\rightarrow\bar{x}$ and for each $k$, $x^k\in \Omega$.  By $x^k\xrightarrow{u}\bar x$, we mean that the sequence $\{x^k\}$ approaches $\bar x$ in direction $u$, i.e., there exist $t_k\downarrow 0, u^k\rightarrow u$ such that $x^k=\bar x+ t_k u^k$. By $f(t)=o(t)$, we mean that $f(t)$ is a function such that $\lim_{t\downarrow0}\frac{f(t)}{t}=0$.
$\mathbb B$ denotes the unit open ball and $\mathbb B_\sigma (\bar x)$ denotes the open ball centered at $\bar x$ with radius equal to $\sigma$. For any  $x\in \mathbb R^n$, we denote by $ \|x\|$ the Euclidean norm.
 For a set $\Omega$, we denote by ${\rm cl}\Omega$, ${\rm int}\Omega$, ${\rm co}\Omega$,  $\Omega^\circ$ and bd$\Omega$ its closure, its interior, its convex hull,  its polar and its boundary, respectively.  By  ${\rm dist}(x,\Omega):=\inf\{\|x-y\||y\in\Omega\}$, we denote the distance from a point $x$ to set $\Omega$. For a single-valued map $\phi:\mathbb R^n\rightarrow\mathbb R^m$, we denote by $\nabla \phi(x)\in \mathbb{R}^{m\times n}$   the Jacobian matrix of $\phi$  at $x$ and for a function $\varphi:\mathbb R^n\rightarrow\mathbb R$, we denote by $\nabla \varphi(x)$ both the gradient and the Jacobian of $\varphi$  at $x$.

We review various concepts of tangent and normal cones below (see, e.g.,  \cite[Definitions 6.1, 6.3 and 6.25, equation 6(19)]{RW}, \cite[Theorem 3.57]{Aub1} and \cite[Definition 2.54]{BS}).
\begin{defn}[Tangent Cone and Normal Cone] \label{tangentnormal}
	Given a set $\Omega\subseteq\mathbb{R}^n$ and a point $\bar{x}\in \Omega$, 
	the tangent/contingent  cone to $\Omega$ at $\bar{x}$ is defined as
	$$T_\Omega(\bar{x}):=\limsup_{t\downarrow0}\frac{\Omega-\bar x}{t}=\left \{d\in\mathbb{R}^n|\exists t_k\downarrow0, d_k\rightarrow d\ \mbox{ s.t. } \bar{x}+t_kd_k\in \Omega\ \forall k\right \}.$$	
	The regular/Clarke tangent cone to $\Omega$ at $\bar x$ is defined as
	\[
	\widehat{T}_\Omega(\bar x):=\liminf_{x\xrightarrow{\Omega}\bar x,t\downarrow0}\frac{\Omega-x}{t}=\left\{d\in\mathbb R^n| \forall t\downarrow0,\forall x\xrightarrow{\Omega}\bar x,\exists y\xrightarrow{\Omega}\bar x\ s.t.\ \frac{y-x}{t}\rightarrow d\right\}.
	\] 
	The regular normal cone,  the limiting normal cone and the Clarke normal cone to $\Omega$ at $\bar{x}$ can be defined as 
\begin{eqnarray*}\widehat{N}_\Omega(\bar{x})&:=&\left \{ \zeta\in \mathbb{R}^n\bigg| \langle \zeta ,x-\bar{x}\rangle \leq o(\|x-\bar x\|) \quad \forall x\in \Omega \right \},\\
	 N_\Omega(\bar{x})&:=&\limsup_{x\xrightarrow{\Omega}\bar x} \widehat{N}_\Omega(x) =\left \{\zeta\in \mathbb{R}^n\bigg| \exists \ x_k\xrightarrow{\Omega}\bar{x},\ \zeta_k{\rightarrow}\zeta\ \text{such that}\ \zeta_k\in\widehat{N}_\Omega(x_k)\ \forall k\right 
	 \},\\
	N^c_\Omega(\bar{x})&:=& {\rm clco} N_\Omega(\bar{x}),  \end{eqnarray*}
	respectively.
	\end{defn}

According to \cite[Formula  (2.88)]{BS}, the  contingent cone can be equivalently written in the following form
\begin{eqnarray}
T_\Omega(\bar{x})&=&\left \{d\in\mathbb{R}^n|\exists t_k\downarrow0, \mbox{\rm dist}(\bar x+t_kd,\Omega)=o(t_k) \right\}. \label{tangentc}
\end{eqnarray}
 By \cite[Definition 6.25 and Theorem 6.26]{RW}, 
the regular tangent cone can be defined as the following equivalent form:
$$\widehat{T}_\Omega(\bar x):=\liminf_{x \rightarrow \bar x,  x\in\Omega}T_\Omega( x).$$

For a point $\bar x\not \in \Omega$, it is conventional to define the corresponding tangent and normal cones to be the emptyset. From definition, we have the following relationships between various tangent cones and normal cones.
	\begin{eqnarray*}
	\widehat{T}_\Omega(\bar{x}) \subseteq 
	{T}_\Omega(\bar{x}),
	&& \widehat{N}_\Omega(\bar{x})  \subseteq   {N}_\Omega(\bar{x})\subseteq   {N}^c_\Omega(\bar{x}).
	\end{eqnarray*}

{
It is well-known that the Clarke normal cone and the regular tangent cone are polar to each other. 
\begin{prop}[Tangent-Normal Polarity]\cite[Theorem 6.28]{RW} For a closed set $\Omega$ and $ \bar x\in \Omega$, one has 
	$$\widehat{T}_\Omega(\bar x)=N^c_\Omega(\bar x)^\circ,\quad  \widehat{T}_\Omega(\bar x)^\circ =N^c_\Omega(\bar x).$$
	In particular, the set $\widehat{T}_\Omega(\bar x)$ is closed and convex.\label{polarity}
\end{prop}}


	The following directional version of limiting normal cone  introduced in \cite[Definition 2.3]{GM} is comprehensively studied by Gfrerer (see e.g., \cite[Definition 2]{Gfr13} and \cite{BHA}).
	\begin{defn}[Directional Normal Cone] 
    Given a set $\Omega \subseteq \mathbb{R}^n$, a point $\bar x \in \Omega$ and a direction $d\in\mathbb{R}^n$, the limiting normal cone to $\Omega$ at $\bar{x}$ in direction $d$ is defined by
    $$N_\Omega(\bar{x};d):=\limsup_{x\xrightarrow{d}\bar x} \widehat{N}_\Omega(x)=\left \{\zeta \in \mathbb{R}^n\bigg| \exists \ t_k\downarrow0, d_k\rightarrow d, \zeta_k\rightarrow\zeta  \mbox{ s.t. } \zeta_k\in \widehat{N}_\Omega(\bar{x}+t_kd_k)\ \forall k \right \},$$ and the Clarke normal cone to $\Omega$ at $\bar{x}$ in direction $d$ is defined by 
    $$N^c_\Omega(\bar{x};d):={\rm clco}N_\Omega(\bar{x};d).$$
\end{defn}
It is obvious that $N^c_{\Omega}(\bar x; 0)=N^c_{\Omega}(\bar x)$, $N^c_{\Omega}(\bar x; d)=\emptyset$ if $d \not \in T_\Omega(\bar x)$ and $N^c_{\Omega}(\bar x;d)\subseteq N^c_\Omega(\bar x)$. It is also obvious that for all $d\in T_\Omega(\bar x) \setminus T_{bd \Omega}(\bar x)$, one has $N^c_\Omega(\bar x;d)=\{0\}$. 
 By \cite[Lemma 2.1]{Gfr14}, the directional and the classical normal cone have the following relationship 
\begin{equation*}
N_\Omega(\bar x;d)\subseteq N_\Omega(\bar x)\cap \{d\}^\perp  \qquad \forall d\in T_\Omega(\bar x), 
\end{equation*} if  $\Omega$ is the union of finitely many convex polyhedral sets,
and the equality holds if $\Omega$ is convex.
{Moreover it is easy to see that $N_\Omega(\bar x;d)$ is homogeneous with degree zero in $d$, i.e.,  $N_\Omega(\bar x;\alpha d)=N_\Omega(\bar x; d)$ for any  $\alpha>0$.}
The following fact proved in  \cite[Proposition 2]{HYZ} is the directional counterpart of outer semicontinuous property of the limiting normal cone mapping (\cite[Proposition 6.6]{RW}).
\begin{prop}
	Given a set $\Omega \subseteq \mathbb{R}^n$, a point $\bar x \in \Omega$ and a direction $d\in\mathbb{R}^n$, one has
$$N_\Omega(\bar{x};d)=\limsup_{x\xrightarrow{d}\bar x} {N}_\Omega(x)=\left \{\zeta \in \mathbb{R}^n\bigg| \exists \ t_k\downarrow0, d_k\rightarrow d, \zeta_k\rightarrow\zeta  \mbox{ s.t. } \zeta_k\in {N}_\Omega(\bar{x}+t_kd_k)\ \forall k \right \}.$$
\end{prop}


Similarly, many classical concepts in variational analysis can have their directional versions. To do this, the following concept of a directional neighborhood is needed.
\begin{defn}[Directional Neighborhood](\cite[formula (7)]{Gfr13}). 
  Given a direction $d\in \mathbb{R}^n$, and positive scalars $\varepsilon,\delta$, the directional neighborhood of direction $d$ is a set defined by 
\begin{equation*}
{\cal V}_{\varepsilon,\delta}(d):=\left \{ \begin{array}{ll}
\{0\} \cup \left \{ z \in \varepsilon\mathbb{B}\setminus \{0\} \mid |\frac{z}{\|z\|} -\frac{d}{\|d\|} |\leq \delta \right \} & \mbox{ if } d\not =0,\\
\varepsilon\mathbb{B} &  \mbox{ if } d =0. \end{array} \right . 
\end{equation*}
\end{defn}

With the directional neighborhood at hand, many classical concepts have been exteneded to their directional versions, and a directional version of the variational analysis has been stimulated. In the sequel, we list some useful results from directional variational analysis.
\begin{defn}[Directional  semicontinuity and continuity] Let $\varphi:\mathbb{R}^n \rightarrow [-\infty,\infty]$ be finite at $\bar x$. We say $\varphi$ is lower semicontinuous (l.s.c.) at $\bar x$ in direction $u$ if 
$$\varphi(\bar x) \leq \liminf_{x\xrightarrow{u}\bar x}\varphi(x). $$ 
We say $\varphi$ is upper semicontinuous (u.s.c.) at $\bar x$ in direction $u$ if 
$$\varphi(\bar x) \geq \limsup_{x\xrightarrow{u}\bar x}\varphi(x). $$ 
We say that $\varphi$ is continuous at $\bar x$ in direction $u$ if 
$$\varphi(\bar x) = \lim_{x\xrightarrow{u}\bar x}\varphi(x). $$
\end{defn}

In this paper we study the directional derivatives in the following sense.
\begin{defn}[Directional derivatives]
Let $\varphi:\mathbb R^n\rightarrow\mathbb{R} $. 
The upper Dini directional derivative of $\varphi$ at $x$ in direction $u$ in Hadamard sense  is
$$
	\varphi'_+(x;u):=\limsup_{t\downarrow 0, u' \rightarrow u}\frac{\varphi(x+tu')-\varphi(x)}{t}.
$$
The lower Dini directional derivative of $\varphi$ at $x$ in direction $u$ in Hadamard sense  is
$$
	\varphi'_-(x;u):=\liminf_{t\downarrow 0, u' \rightarrow u}\frac{\varphi(x+tu')-\varphi(x)}{t}.
$$
If both limits are finite and equal, we have then the directional derivative of $\varphi$ at $x$ in direction $u$ in Hadamard sense 
$$
	\varphi'(x;u):=\lim_{t\downarrow 0, u' \rightarrow u}\frac{\varphi(x+tu')-\varphi(x)}{t}.
$$
 We say that $\varphi$ is Hadamard directionally differentiable at $x$ if its directional derivative  in Hadamard sense exists at any direction.
\end{defn}
It is clear that if a function is Hadamard directionally differentiable, then it is directionally differentiable in the usual sense as well, i.e.,
$$
	\varphi'(x;u)=\lim_{t\downarrow 0}\frac{\varphi(x+tu)-\varphi(x)}{t}.
$$ So the concept of Hadamard directional differentiability is stronger than the usual directional differentiability.

Robinson's constraint qualification  (Robinson's CQ)
 (see e.g., \cite{BS})
 is an important constraint qualification and was originally defined for the set-constrained problem $(P_{\bar x})$ where $C$ is a convex set and the constraint region is convex.   For a general set-constrained problem $(P_{\bar x})$, we say that Robinson's CQ holds at $ y\in\mathcal F(\bar x)$ if  \begin{equation} \label{RCQ} \nabla_y P(\bar x, y)^T\lambda=0,\ \lambda\in N^c_C(P(\bar x, y))\implies\lambda=0.\end{equation}
 Let $\bar y \in {\cal F}(\bar x)$. We denote the set of  Clarke multipliers and the limiting  multipliers  at $\bar y$ for problem $(P_{\bar x})$ respectively by  
\begin{eqnarray*}
\Lambda^c(\bar x,\bar y)&:=&\left \{\lambda\in N^c_C(P(\bar x,\bar y))|\nabla_y f(\bar x,\bar y)+\nabla_y P(\bar x,\bar y)^T\lambda=0 \right \},\\
\Lambda(\bar x,\bar y)&:=&\left \{\lambda\in N_C(P(\bar x,\bar y))|\nabla_y f(\bar x,\bar y)+\nabla_y P(\bar x,\bar y)^T\lambda=0 \right \}.\end{eqnarray*}
  It is known that (\ref{RCQ}) holds if and only if the set of Clarke multipliers $\Lambda^c(\bar x, y)$ is nonempty and compact. By the tangent-normal polarity in Proposition \ref{polarity}, it is easy to see that  Robinson's CQ holds at $y$ if and only if 	\begin{align}\label{Robinsonnew}
			\nabla_y P(\bar x, y)\mathbb R^m+{\widehat{T}_C(P(\bar x, y))}=\mathbb R^p.
		\end{align}

Consider a conic linear program in the form of \cite[(2.345)]{BS}:
\begin{align}
\label{primalP}
(P)~~~~~~~\min_{x} &  \ \ \langle \alpha, x \rangle+c  \nonumber \\
s.t. &\ \  Ax +b \in K,
\end{align}
where $A\in\mathbb R^n\times\mathbb R^m$ is a matrix, $\alpha, b\in\mathbb R^m$ are vectors, $c$ is a scalar and $K\subseteq\mathbb R^m$ is a closed convex cone. According to \cite[(2.347)]{BS}, its dual problem can be written in the form
\begin{align}
\label{dualD}
(D)~~~~~~~\max_{\lambda \in K^\circ} & \ \  \langle \lambda, b \rangle +c \nonumber \\
s.t. &\ \   A^T\lambda +\alpha =0.
\end{align}
Define the value function for the perturbed problem of (P) as
$$v(y):=\inf_x \{ \langle \alpha, x \rangle+c: Ax+b+y\in K\}.$$
It is known that the function $v(y)$ is convex.
By \cite[Definition 2.146]{BS}, we say that problem (P) is calm if $v(0)$ is finite and the convex function  $v(y)$ is subdifferentiable at $y=0$, i.e., $\partial v(0)\not =\emptyset.$
By \cite[Proposition 2.147]{BS}, if (P) is calm then there is no duality gap and the set of optimal solutions of the dual problem (D) is nonempty, i.e.,
$$\inf_x \{ \langle \alpha, x \rangle+c: Ax+b\in K\}=\max_{\lambda\in K^\circ} \{ \langle \lambda, b \rangle +c : A^T\lambda +\alpha=0\}.$$

{It is known that if $v(0)$ is finite and  $K$ is a convex polyhedral set then (P) is calm (see \cite{Robinson81} and \cite[Proposition 2.186]{BS}) and consequently the above strong duality holds automatically.} Another sufficient condition for the calmness of (P)  is    the regularity condition
\begin{equation*}
  0\in {\rm int} \{ A (\mathbb{R}^n)+b-K\} .\end{equation*}
The above condition is equivalent to that Robinson's CQ holds at each feasible solution 
 of problem (P), i.e., for any feasible point $x_0$ of (P),
$$0\in {\rm int} \{Ax_0+ A (\mathbb{R}^n-x_0)+b-K\}$$
 which is equivalent to
$$ A(\mathbb{R}^n) -T_K(Ax_0+b)  =\mathbb{R}^m$$
or equivalently
$$A^T\lambda =0, \quad \lambda \in N_K(Ax_0+b) \Longrightarrow \lambda=0. $$

{By \cite[Theorem 2.165]{BS} and the comments in \cite[Page 113]{BS}, one has the following duality theorem for the conic linear program under Robinson's CQ.
\begin{thm}\label{dt}
	Suppose Robinson's CQ holds at each feasible solution of  (P). Then there is no duality gap between problems (P) and (D), i.e., val(P)=val(D). Moreover, if the optimal value of (P) is finite, then the optimal solution set of the dual problem (D) is a nonempty, convex and compact set.
\end{thm}}

\section{Directional Robinson stability}

In this section, we study a concept of stability subject to a directional  perturbation.  
Although the result is mainly for the study of directional derivative of the value function in the next section, it  is of independent interest. 

For fixed $(\bar x,\bar y) \in {\rm gph} {\cal F}:=\{(x,y)\mid y\in {\cal F}(x)\}=\{(x,y)\mid P(x,y)\in C\}$, we say that the set-valued map $y\rightrightarrows P(\bar x,y)-C$ or the system $P(\bar x,y)\in C$ is metrically subregular at $\bar y$ or the local error bound holds at $\bar y$ if there exist constants $\kappa_{\bar x} \geq 0, \varepsilon\geq 0$ such that 
$${\rm dist}(y, {\cal F}(\bar x))\leq \kappa_{\bar x} {\rm dist}(P(\bar x,y),C) \quad \forall y\in \bar y+  \varepsilon \mathbb{B}. 
$$ 
Note that the modulus of the metric subregularity  depends on the parameter $x$. In stability analysis of a system with a perturbation, it  is important to study when this modulus is independent of the parameter locally. 
As early as in 1976, Robinson \cite{Robinson} studied this property for the case where $C$ is a convex cone under the name of ``stability'' and proved that  Robinson's CQ  for problem $(P_{\bar x})$ is a sufficient condition for this kind of stability. Recently Gfrerer and Mordukhovich \cite{HM} extended this concept to the general system and call it the Robinson stability. 
 Usually we are interested in stability with a perturbation from a point $\bar x$ lying in a neighborhood of $\bar x$. However sometimes we only  concern about a perturbation from $\bar x$  which lies in a directional neighborhood of $\bar x$.  
 In particular when we study the directional derivative of the value function we are only interested in a perturbation from $\bar x$ in the given direction. Based on this consideration, 
in \cite[Definition 4.7]{BY}, a directional version of the Robinson stability property (see e.g., \cite[Definition 1.1]{HM}) was introduced for the system of inequalities. 
We now define it for our more general system in the same way.
\begin{defn}[Directional Robinson stability]
We say that ${\cal F}$ satisfies the Robinson stability ${\rm (RS)}$  at $(\bar x,\bar y)\in {\rm gph }{\cal F}$ in a direction $u \in\mathbb R^n $ if there exist {scalars $\kappa\geq0$ and $\varepsilon,\delta>0$} such that 
$${\rm dist}(y, {\cal F}(x))\leq \kappa {\rm dist}(P(x,y),C) \quad \forall (x,y) \in (\bar x,\bar y)+ {\mathcal V}_{\varepsilon,\delta}(u)\times \varepsilon \mathbb{B}.
$$ When $u=0$ in the above definition, we say that ${\cal F}$ satisfies  RS  at $(\bar x,\bar y)$.
\end{defn}
RS in direction $u$ means that the   local  error bound condition holds at $\bar y$  uniformly in a  neighborhood of $\bar x$ in direction $u$. Hence RS at $(\bar x,\bar y)\in {\rm gph }{\cal F}$ in direction $u$ implies that that the system $P(\bar x, y)\in C$ is metrically subregular at $\bar y$.

Since $P$ is smooth, it is known that the mapping $y \rightrightarrows P(\bar x,y)-C$ is metrically regular at $\bar y$ if and only if the No Nonzero Abnormal Multiplier constraint qualification (NNAMCQ) for problem $(P_{\bar x})$ holds at $\bar y$:
\begin{equation} \label{NNAMCQ} \nabla_y P( \bar x,\bar  y)^T\lambda=0,\ \lambda\in {N}_C(P( \bar x, \bar  y))\implies\lambda=0.\end{equation} Note that when $C$ is closed and convex, the NNAMCQ is equivalent to Robinson's CQ defined in (\ref{RCQ}) and  for the case of nonlinear program, NNAMCQ is equivalent to MFCQ for nonlinear programs.

 In the following proposition,  Gfrerer and Mordukhovich  \cite[Corollary 3.7]{HM} extended Robinson's result that Robinson's CQ implies RS to the more general system where $C$ is a closed set.
\begin{prop}\label{RCQRS}
\cite[Corollary 3.7]{HM}
 Let $\bar y\in {\cal F}(\bar x)$. If the NNAMCQ for problem $(P_{\bar x})$ holds at $\bar y$, then RS holds at $(\bar x,\bar y)$.
 \end{prop}
 
 By the well-known Ronbinson's polyhedral multifunction theory \cite{Robinson81}, we know that a polyhedral multifunction is upper Lipschitz continuous. It turns out that  a parametric polyhedral multifunction is Robinson stable if the parameter $x$ is separated from the variable $y$.
 \begin{prop}\label{prop3.2} Let $P(x,y)=a(x)+By+c$ with $a:\mathbb R^n\rightarrow\mathbb R^p$ being continuous, $B\in\mathbb R^{p\times m}, c\in \mathbb{R}$, $C$ is the union of finitely many convex polyhedral sets and $u$ be a direction. If $\mathcal F(x)$ is nonempty {for all $x \in {\mathcal V}_{\varepsilon,\delta}(u)$ for some $\varepsilon, \delta>0$}, then ${\cal F}$ satisfies  RS  at each $(\bar x,\bar y)\in {\rm gph} {\cal F}$ in direction $u$.
\end{prop}
\beginproof
By \cite[Proposition 1]{Robinson81}, for {each $x$} if ${\cal F}(x)\not = \emptyset$, then the metric subregularity holds, i.e.,
\begin{equation}\label{MS}
{\rm dist}(y,\mathcal F(x))\leq\kappa{\rm dist}(P(x,y),C) \qquad \forall y\in \mathbb{R}^m.
\end{equation}Furthermore, the modulus $\kappa$ depends only on the matrix $B$ and the set $C$ (see \cite[Proposition 1]{Robinson81} and \cite[proof of Theorem 1]{Walkup69} for details). Hence if  ${\cal F}(x)\not = \emptyset$ { for any $x$ in ${\mathcal V}_{\varepsilon,\delta}(u)$, (\ref{MS}) holds for all $x \in {\mathcal V}_{\varepsilon,\delta}(u)$, which  means that RS holds in direction $u$.
\endproof

{We say that the first order sufficient condition for metric subregularity (FOSCMS) of the system $P(\bar x,y)\in C$ at $\bar y$  holds in direction $\nabla_y P(\bar x,\bar y)v $  if 
\begin{equation*}
	\nabla_y P(\bar x,\bar y)^T\lambda=0,\ \lambda\in N_C(P(\bar x,\bar y); \nabla_y P(\bar x,\bar y)v  )\implies\lambda=0. 
	\end{equation*} By Gfrerer and Klatte \cite[Corollary 1]{GKlatte16} if FOSCMS of the system $P(\bar x,y)\in C$ at $\bar y$  holds in each direction $\nabla_y P(\bar x,\bar y)v \in T_C(P(\bar x,\bar y))$ with $v\not =0$, then the set-valued map $y\rightrightarrows P(\bar x,y)-C$ is metrically subregular  at $\bar y$. In fact, the set of all $v$ satisfying $$\nabla_y P(\bar x,\bar y)v \in T_C(P(\bar x,\bar y))$$  is the set of linearization cone for the feasible region of problem $(P_{\bar x})$ at $\bar y$. FOSCMS holding in all nonzero directions in the linearization cone is weaker than the NNAMCQ (\ref{NNAMCQ}) and it does not imply RS. It is natural to ask under what extra conditions, FOSCMS implies RS. This question is addressed by Gfrerer and Mordukhovich in 
	\cite{HM}. Motivated by their research, we now extend their results to the directional RS.}

The reason that FOSCMS in all nonzero directions in the linearization cone does not imply RS is that the dependence of the parameter $x$ is not reflected in those directions. Let $u$ be  given. In order to modify FOSCMS so that it implies RS, we  consider the following direction.
Define the image directional derivative of $P$ with respect to $x$ at $(\bar x,\bar y)$ in direction $u$ as the closed cone $${\rm Im}D_xP(\bar x,\bar y;u):=\left\{\alpha v\left|\alpha\geq0, \exists t_k\downarrow0, u^k\rightarrow u\ s.t.\ v=\lim_{k\rightarrow\infty}\frac{P(\bar x+t_ku^k,\bar y)-P(\bar x,\bar y)}{\|P(\bar x+t_ku^k,\bar y)-P(\bar x,\bar y)\|}\right.\right\}.$$
Obviously the convex cone generated by the vector  $\nabla_xP(\bar x,\bar y)u$ is a subset of the set ${\rm Im}D_xP(\bar x,\bar y;u)$,  i.e., $$\{\alpha\nabla_xP(\bar x,\bar y)u|\alpha\geq0\}\subseteq {\rm Im}D_xP(\bar x,\bar y;u).$$ The equality holds if $\nabla_xP(\bar x,\bar y)u\neq0$.
 In general, the strict inclusion may hold. {When $u=0$, ${\rm Im}D_xP(\bar x,\bar y):={\rm Im}D_xP(\bar x,\bar y;0)$ coincides with the image derivative introduced in Gfrerer and Mordukhovich \cite[formula (3.8)]{HM} with $\zeta(x):=\|x-\bar x\|$. 
 
For $(x,y) \in {\rm gph} {\cal F}$ and given direction $u$, define the closed cone
	\begin{eqnarray}
&&{D(x,y;u):=}
	 \left \{d+\nabla_y P( x, y)v \left|\begin{array}{l}
	0\not =(d,v)\in ({\rm Im}D_xP( x, y;u)\times \mathbb{R}^m),\\
	 d+\nabla_y P( x, y)v \in 
	T_C(P(  x, y))\end{array} \right. \right \}.
	 \label{eqn12}
	 \end{eqnarray}
	By taking $d=0$ in the above, it is obvious that the cone $D(\bar x,\bar y;u)$ contains  all directions $\nabla_y P(\bar x,\bar y)v \in T_C(P(\bar x,\bar y))$ with $v\not =0$.

We now present  sufficient conditions for the directional RS using the first order information. The result extends Gfrerer and Mordukhovich \cite[Corollary 3.6]{HM} in that the image directional derivative ${\rm Im}DP_x(\bar x,\bar y;u)$ instead of the image derivative ${\rm Im}D_xP(\bar x,\bar y)$ is used. Hence in the case where $u=0$ our result recovers  \cite[Corollary 3.6]{HM} applied to the parametric program $(P_x)$}. Since for $u\not =0$, ${\rm Im}D_xP(\bar x,\bar y;u)$ is in general strictly contained in ${\rm Im}D_xP(\bar x,\bar y)$, the sufficient condition for the directional RS in Proposition \ref{SRS} is weaker than the one proposed in \cite[Corollary 3.6]{HM}. Naturally the conclusion here is only the directional RS while the conclusion in \cite[Corollary 3.6]{HM} is the full/nondirectional RS. Using  the proof technique  of \cite[Theorem 3.5]{HM}, we can now prove the following result.
{Since a direction $\nabla_y P(\bar x,\bar y)v \in T_C(P(\bar x,\bar y))$ with $v\not =0$ lies in $D(\bar x,\bar y;u)$ with $d=0$, it is obvious FOSCMS  (\ref{eqn(15)}) holding in all directions $w\in D(\bar x,\bar y;u)$ is stronger than FOSCMS holding in all directions $\nabla_y P(\bar x,\bar y)v \in T_C(P(\bar x,\bar y))$ with $v\not =0$ in that the perturbation in $x$ is now considered.}
\begin{prop}[First-order sufficient conditions for the directional RS]\label{SRS}Let $\bar y\in \mathcal F(\bar x)$
 and $u\in \mathbb{R}^n$ be a given direction. Suppose that for  every $d\in{\rm Im}D_xP(\bar x,\bar y;u)$ and every $\tau_k\downarrow 0$,  there exists $v$ such that 
	\begin{equation}\label{geom}
	\liminf_{k\rightarrow \infty }
	\frac{{\rm dist}(P(\bar x,\bar y)+\tau_k(d+\nabla_yP(\bar x,\bar y)v),C)}{\tau_k}
	=0,
	\end{equation}
	and {FOSCMS of the system $P(\bar x,y)\in C$ at $\bar y$ holds in every direction $w\in D(\bar x,\bar y;u)$. i.e,}
	\begin{equation}
	\nabla_y P(\bar x,\bar y)^T\lambda=0,\ \lambda\in N_C(P(\bar x,\bar y);w )\implies\lambda=0. \label{eqn(15)}
	\end{equation}
  Then  ${\cal F}$ satisfies RS  at $(\bar x,\bar y)$ in direction $u$.


\end{prop}
\beginproof
 We prove the result by contradiction. Assume  RS fails at $(\bar x,\bar y)$ in direction $u$. Then there exist sequences $t_k\downarrow0, u^k\rightarrow u, y^k\rightarrow\bar y$ such that 
\begin{equation}\label{assu}
{\rm dist}(y^k,\mathcal F(\bar x+t_ku^k))>k{\rm dist}(P(\bar x+t_ku^k,y^k),C).
\end{equation}
Without loss of generality we may assume that for any $k\in\mathbb N$, there exists $\epsilon_k>0$ such that 
$ \|\nabla_yP(x,y)-\nabla_yP(\bar x,\bar y)\|\leq1/k$ for any $(x,y)\in\mathbb B_{\epsilon_k}(\bar x,\bar y)$, and
  $$\|\bar x+t_ku^k-\bar x\|<\epsilon_k/2, \|y^k-\bar y\|<\min\{1/k,\epsilon_k/2\}, \|P(\bar x+t_ku^k,y^k)-P(\bar x,\bar y)\|\leq1/k^2.$$
Define 
\begin{equation}\label{sigma} \sigma_k:=1/(k^2{\rm dist}(P(\bar x+t_ku^k,y^k),C)).\end{equation} Let $(\bar y^k,\bar\xi^k)
$ be an optimal solution to the minimization problem
\begin{equation}\label{proj}
\min_{y,\xi}	\varphi_k(y,\xi):=\|\xi\|+\sigma_k\|y-y^k\|^2\quad s.t.\ P(\bar x+t_ku^k,y)+\xi\in C.
\end{equation} We claim that $\bar \xi^k\neq0$. To the contrary, suppose that $\bar \xi^k=0$.
 Take $\xi^k$ such that $P(\bar x+t_ku^k,y^k)+\xi^k\in C$ and 
$
\|\xi^k\|={\rm dist}(P(\bar x+t_ku^k,y^k),C).$
\begin{eqnarray*}
{\rm dist}(y^k,\mathcal F(\bar x+t_ku^k))^2 &\leq &\| y^k-\bar y^k\|^2 \mbox{ since } \bar y^k \in \mathcal F(\bar x+t_ku^k)\\
&\leq & \frac{1}{\sigma_k}{\rm dist}(P(\bar x+t_ku^k,y^k),C) \quad \mbox{since } 
\varphi_k(\bar y^k,\bar\xi^k)\leq\varphi_k(y^k,\xi^k)\\
&=& k^2{\rm dist}(P(\bar x+t_ku^k,y^k),C)^2 \quad \mbox{ by } (\ref{sigma})
\end{eqnarray*} contradicting (\ref{assu}). The contradiction proves that $\bar \xi^k\neq0$. Furthermore, since
\[
\sigma_k\|\bar y^k-y^k\|^2\leq \|\bar \xi^k\|+\sigma_k\|\bar y^k-y^k\|^2= \varphi_k(\bar y^k,\bar\xi^k)\leq\varphi_k(y^k,\xi^k)={\rm dist}(P(\bar x+t_ku^k,y^k),C),
\]
and
\begin{eqnarray}
\|\bar y^k-y^k\| & \leq & k{\rm dist}(P(\bar x+t_ku^k,y^k),C) \label{converg}\\
&\leq & k\|P(\bar x+t_ku^k,y^k)-P(\bar x,\bar y)\|\leq\frac{1}{k} \nonumber,
\end{eqnarray} we have $\|\bar y^k-\bar y\|\leq\|\bar y^k-y^k\|+\|y^k-\bar y\|\leq2/k$ and $P(\bar x+t_ku^k,\bar y^k)\rightarrow P(\bar x,\bar y)$ as $k\rightarrow\infty$.

Next we consider the direction along which the sequence $(\bar x+t_ku^k,\bar y^k)$ converges to $(\bar x,\bar y)$.
Define $\tau_k:=\|\bar y^k-\bar y\|+\|P(\bar x+t_ku^k,\bar y)-P(\bar x,\bar y)\|$ for $k\in\mathbb N$. Then $\tau_k\downarrow0$ as $k\rightarrow\infty$ and without loss of generality there exists a vector $(\bar v,\bar d)\in\mathbb R^m\times\mathbb R^p$ such that
\begin{equation}\label{dir}
\bar v:=\lim_k\frac{\bar y^k-\bar y}{\tau_k}\ \mbox{and}\ \bar d:=\lim_k\frac{P(\bar x+t_ku^k,\bar y)-P(\bar x,\bar y)}{\tau_k}.
\end{equation}
One can easily obtain $(\bar v,\bar d)\neq(0,0)$ with $\bar d\in{\rm Im}D_xP(\bar x,\bar y;u)$. To show that $\bar d+\nabla_yP(\bar x,\bar y)\bar v\in T_C(P(\bar x,\bar y))$, {by (\ref{tangentc}) it suffices to  prove that } 
\begin{equation} {\rm dist}(P(\bar x,\bar y)+\tau_k(\bar d+\nabla_yP(\bar x,\bar y)\bar v),C)=o(\tau_k).\label{equivdef} \end{equation}

By condition (\ref{geom}), for $\bar d$ and the sequence $\{\tau_k\}$ there exists $v\in\mathbb R^m$ such that
\[
\liminf_{k\rightarrow \infty} \frac{{\rm dist}(P(\bar x,\bar y)+\tau_k(\bar d+\nabla_yP(\bar x,\bar y)v),C)}{\tau_k}=0.
\]
 Since
\begin{align}\label{diff}
	&\lim_k\frac{\|P(\bar x+t_ku^k,\bar y+\tau_kv)-(P(\bar x,\bar y)+\tau_k(\bar d+\nabla_yP(\bar x,\bar y)v))\|}{\tau_k}\notag\\
	=&\lim_k\frac{\|P(\bar x+t_ku^k,\bar y+\tau_kv)-(P(\bar x+t_ku^k,\bar y)+\tau_k\nabla_yP(\bar x,\bar y)v)\|}{\tau_k}\notag\\
	=&\lim_k\frac{\|\int_0^1[\nabla_yP(\bar x+t_ku^k,\bar y+s(\bar y+\tau_k v-\bar y))-\nabla_yP(\bar x,\bar y)]\tau_kv{\rm d}s\|}{\tau_k}\notag\\
	\leq&\lim_k\frac{2}{k}\|v\|=0,
\end{align} where the first equality follows from (\ref{dir}), one has
\[
\liminf_{k\rightarrow \infty} \frac{{\rm dist}(P(\bar x+t_ku^k,\bar y+\tau_kv),C)}{\tau_k}=0.
\]
Let $\xi^k$ be such that $P(\bar x+t_ku^k,\bar y+\tau_kv)+\xi^k\in C$ and $\|\xi^k\|={\rm dist}(P(\bar x+t_ku^k,\bar y+\tau_kv),C)$. Then $\|\xi^k\|=o(\tau_k)$.

Since $\varphi(\bar y^k,\bar \xi^k)\leq\varphi(\bar y+\tau_kv,\xi^k)$ 
one has
\begin{align}\label{esti}
	\|\bar\xi^k\|\leq&  \|\xi^k\|-\sigma_k\|\bar y^k-y^k\|^2+ \sigma_k\|\bar y+\tau_kv- y^k\|^2\notag\\
	\leq& \|\xi^k\|+\sigma_k(2\langle\bar y^k- y^k,\bar y+\tau_kv-\bar y^k\rangle+\|\bar y+\tau_kv-\bar y^k\|^2)\notag\\
	\leq&\|\xi^k\|+2\sigma_k\|\bar y^k- y^k\|\cdot\|\bar y+\tau_kv-\bar y^k\|+\sigma_k\|\bar y+\tau_kv-\bar y^k\|^2.
\end{align}
By (\ref{converg}) we have
\begin{equation}\label{esti2}
	\sigma_k\|\bar y^k- y^k\|\leq\frac{1}{k^2{\rm dist}(P(\bar x+t_ku^k,y^k),C)}k{\rm dist}(P(\bar x+t_ku^k,y^k),C)=\frac{1}{k}\rightarrow0\ \mbox{as}\ k\rightarrow\infty.
\end{equation}
Since $\|\bar y+\tau_kv-\bar y^k\|\leq\|\bar y+\tau_kv-\bar y\|+\|\bar y^k-\bar y\|\leq \tau_k(\|v\|+2\|\bar v\|)$ and combining (\ref{esti}), (\ref{esti2})
and the facts $\|\bar \xi^k\|\leq{\rm dist}(P(\bar x+t_ku^k,y^k),C)$, $\|\xi^k\|=o(\tau_k)$ yield
\begin{align}\label{clo}
	&\lim_k\frac{\|\bar \xi^k\|}{\tau_k}\notag\\
	\leq&\lim_k\frac{\min\{\|\xi^k\|+2\sigma_k\|\bar y^k- y^k\|\cdot\|\bar y+\tau_kv-\bar y^k\|+\sigma_k\|\bar y+\tau_kv-\bar y^k\|^2,{\rm dist}(P(\bar x+t_ku^k,y^k),C)\}}{\tau_k}\notag\\
	\leq&\lim_k\frac{\|\xi^k\|}{\tau_k}+\lim_k\frac{2\sigma_k\|\bar y^k- y^k\|\cdot\|\bar y+\tau_kv-\bar y^k\|}{\tau_k}\notag\\
	& +\lim_k\frac{\min\{\sigma_k\|\bar y+\tau_kv-\bar y^k\|^2,{\rm dist}(P(\bar x+t_ku^k,y^k),C)\}}{\tau_k}\notag\\
	\leq&0+\lim_k\frac{2/k\tau_k(\|v\|+2\|\bar v\|)}{\tau_k}+\lim_k\frac{\sqrt{\sigma_k{\rm dist}(P(\bar x+t_ku^k,y^k),C)}\|\bar y+\tau_kv-\bar y^k\|}{\tau_k}\notag\\
	\leq&\lim_k\frac{2}{k}(\|v\|+2\|\bar v\|)=0,
\end{align} where the last inequality follows from (\ref{sigma}).

Combining (\ref{clo}) and (\ref{diff}), one can easily obtain
\begin{align*}
&{\rm dist}(P(\bar x,\bar y)+\tau_k(\bar d+\nabla_yP(\bar x,\bar y)\bar v),C)\\
\leq&\|P(\bar x,\bar y)+\tau_k(\bar d+\nabla_yP(\bar x,\bar y)\bar v)-P(\bar x+t_ku^k,\bar y+\tau_k\bar v)\|+{\rm dist}(P(\bar x+t_ku^k,\bar y+\tau_k\bar v),C)\\
\leq&o(\tau_k)+{\rm dist}(P(\bar x+t_ku^k,\bar y^k),C)
=o(\tau_k)+\|\bar\xi^k\|=o(\tau_k),
\end{align*} where {the second inequality follows from (\ref{diff}) and (\ref{dir})} and the last equality is obtained in (\ref{clo}). By (\ref{equivdef}), $\bar d+\nabla_yP(\bar x,\bar y)\bar v\in T_C(P(\bar x,\bar y))$.

Finally, we are ready to deduce a contradiction. The constraint system for problem (\ref{proj}) can be written as  $H_k(y,\xi)\in C$ where $H_k(y,\xi):=P(\bar x+t_ku^k,y)+\xi$. Since the Jacobian $\nabla H_k(y,\xi)$ has full row rank, the   necessary optimality conditions for (\ref{proj}) holds at $(\bar y^k,\bar\xi^k)$, i.e., there exists a multiplier $\lambda_k\in N_C(P(\bar x+t_ku^k,\bar y^k)+\bar\xi^k)$ such that
\begin{equation}\label{cont}
2\sigma_k(\bar y^k-y^k)+\nabla_yP(\bar x+t_ku^k,\bar y^k)^T\lambda^k=0\ \mbox{and}\ \frac{\bar\xi^k}{\|\bar \xi^k\|}+\lambda^k=0.
\end{equation}
Since $\|\lambda^k\|=1$, passing to a subsequence if necessary, there exists $0\neq\bar \lambda=\lim_k\lambda^k$. Then taking the limit of (\ref{cont}) as $k\rightarrow\infty$, one has
\[
\nabla_yP(\bar x,\bar y)^T\bar\lambda=0,\ \bar\lambda\in N_C(P(\bar x,\bar y);\bar d+\nabla_yP(\bar x,\bar y)\bar v)
\]contradicting condition (\ref{eqn(15)}). The proof is complete.
\endproof

Since when $\nabla_xP(\bar x,\bar y)u\neq0$, $\{\alpha\nabla_xP(\bar x,\bar y)u|\alpha\geq0\}={\rm Im}D_xP(\bar x,\bar y;u)$, in this case Proposition \ref{SRS} has the following form.
\begin{cor}\label{cor3.1}Let $\bar y\in \mathcal F(\bar x)$ and
 $u\in \mathbb{R}^n$ be a given direction such that $\nabla_xP(\bar x,\bar y)u\neq0$. Suppose that for every $\alpha\geq0$
 and every $\tau_k\downarrow 0$,  
there exists $v$ such that 
		\begin{equation}\label{geom1}
			\liminf_{k\rightarrow \infty }
			\frac{{\rm dist}(P(\bar x,\bar y)+\tau_k\nabla P(\bar x,\bar y)(\alpha u, v),C)}{\tau_k}
			=0,
		\end{equation}
		and 
		FOSCMS of the system $P(\bar x,y)\in C$ at $\bar y$ holds in direction  $\nabla P(\bar x,\bar y)(\alpha u, v) \in T_C(P(\bar x,\bar y))$ for every $(0,0)\neq(\alpha,v)$ with $\alpha\geq 0$, i.e.,
		\begin{equation*}
			\nabla_y P(\bar x,\bar y)^T\lambda=0,\ \lambda\in N_C(P(\bar x,\bar y);\nabla_yP(\bar x,\bar y)(\alpha u, v))\implies\lambda=0. 
		\end{equation*} 
Then  ${\cal F}$ satisfies  RS  at $(\bar x,\bar y)$ in direction $u$.
\end{cor}

\begin{remark}[Assumption Verification]\label{Remark3.1}

The verification of Condition (\ref{geom}) in Proposition \ref{SRS} may not be easy to  verify. However it holds automatically in several usual cases.
In fact, by the proof of \cite[Corollary 3.7]{HM}, if NNAMCQ for problem $(P_{\bar x})$ holds at $\bar y\in {\cal F}(\bar x)$, then condition (\ref{geom}) in Proposition \ref{SRS} always holds. By \cite[Remark 3.8]{HM}, in case where $C$ is the union of finitely many polyhedral sets, {if the vector $(d,v)$ fulfill the condition 
\begin{equation}\label{assumpv}
	d+\nabla_y P(\bar x,\bar y)v\in T_C(P(\bar x,\bar y)),
\end{equation} then there exists $\bar t>0$ with 
$$
P(\bar x,\bar y)+t(d+\nabla_y P(\bar x,\bar y)v )\in C\ \mbox{for all}\ t\in[0,\bar t].
$$} Consequently, condition (\ref{geom}) in Proposition \ref{SRS} { reduces to the feasibility of system (\ref{assumpv})} when $C$ is the union of finitely many polyhedral sets. Similarly in the case where $\nabla_xP(\bar x,\bar y)u\neq0$,   if either NNAMCQ  for problem $(P_{\bar x})$ holds at $\bar y\in {\cal F}(\bar x)$ or $C$ is the union of finitely many polyhedral sets and the system (\ref{assumpv})  with $d=\alpha \nabla_x P(\bar x,\bar y)u$  is feasible, then  condition (\ref{geom1}) in Corollary \ref{cor3.1} always hold. 
\end{remark}
Hence we have the following sufficient condition for RS.
\begin{cor}Suppose $C$ is the union of finitely many polyhedral sets, $\bar y\in \mathcal F(\bar x)$ and  for every $d\in{\rm Im}D_xP(\bar x,\bar y;u)$, there exists $v$ satisfying
	\[
		d+\nabla_y P(\bar x,\bar y)v\in T_C(P(\bar x,\bar y)).
	\] Suppose FOSCMS  of the system $P(\bar x,y)\in C$ at $\bar y$ holds  i.e., \begin{equation*}
	\nabla_y P(\bar x,\bar y)^T\lambda=0,\ \lambda\in N_C(P(\bar x,\bar y));w) \implies\lambda=0
	\end{equation*} for every $w\in D(x,y;u)$  as defined in (\ref{eqn12}).
	 Then  ${\cal F}$ satisfies RS  at $(\bar x,\bar y)$ in direction $u$. {Moreover if $\nabla_xP(\bar x,\bar y)u\neq0$, then any vector $w$ in $ D(x,y;u)$ can be written as $w=\nabla P(\bar x,\bar y)(\alpha u, v) $ for some $\alpha \geq 0$ and $v \in \mathbb{R}^m$ such that  $(0,0)\neq(\alpha,v)$ and $\nabla P(\bar x,\bar y)(\alpha u, v) \in T_C(P(\bar x,\bar y))$.}
\end{cor}

From the proof of Proposition \ref{SRS}, we can see that the first  order sufficient conditions for RS for splitting systems in \cite[Theorems 3.5]{HM} can be similarly extended with ${\rm Im}D_xP(\bar x,\bar y)$ replaced by  ${\rm Im}D_xP(\bar x,\bar y;u)$. 
Note that second order sufficient conditions for RS for splitting systems  can  also be found in \cite[Theorem 4.3 and Corollary 4.4]{HM}. 

At the end of this section, we give an application of directional RS which will be useful in the next section. In  \cite[Lemma 4.2]{BY}, it was shown that for a parametric nonlinear program, RS in direction {$u$} implies that the linearization system  (\ref{LS})  is feasible. We now extend the result to our more general case.
\begin{prop}\label{LinearizC} Let $\bar y\in {\cal F}(\bar x)$ and $u\in \mathbb{R}^n$. Suppose ${\cal F}$ satisfies RS at $(\bar x,\bar y)$ in direction $u$. Then there exists $v$ satisfying the  system \begin{equation}\label{LS}
	\nabla P(\bar x,\bar y)(u, v) \in T_C(P(\bar x,\bar y)) . 
	\end{equation}
\end{prop}
{\beginproof
Since RS holds at $(\bar x,\bar y)$ in direction $u$, by Definition 3.1, there exist positive scalars $\kappa,\epsilon,\delta$ such that, for any $x\in \bar x+\mathcal V_{\epsilon,\delta}(u)$,
\begin{equation}\label{16}
{\rm dist}(\bar y,\mathcal F(x))\leq\kappa{\rm dist}(P(x,\bar y),C)\leq\kappa\|P(x,\bar y)-P(\bar x,\bar y)\|\leq L_p\kappa\|x-\bar x\|,
\end{equation}where $L_p$ is the Lipschitz modulus of $P(x,\bar y)$ around $\bar x$. Then for any sequences $t_k\downarrow0, u^k\rightarrow u$, by (\ref{16}), we can find a sequence $y^k\in\mathcal F(\bar x+t_ku^k)$ such that $$\|\bar y-y^k\|\leq L_p\kappa\|\bar x+t_ku^k-\bar x\|,$$ which implies that $y^k\rightarrow\bar y$. Since $\{(y^k-\bar y)/t_k\}$ is bounded, taking a subsequence if necessary, we can find $v\in\mathbb R^m$ such that $v^k:=(y^k-\bar y)/t_k$ converges to $v$. Then 
$$(P(\bar x+t_ku^k,y^k)-P(\bar x,\bar y))/t_k \rightarrow  \nabla P(\bar x,\bar y)(u,v).$$ Since $P(\bar x+t_ku^k,y^k)\in C$, it follows that  $\nabla P(\bar x,\bar y)(u,v)\in T_C(P(\bar x,\bar y))$.
\endproof}


\section{Directional derivatives using the {directional solution set}}
The main purpose of this section is to extend the classical results on the directional derivative to the general parametric problem $(P_x)$ where the solution set is replaced by a directional solution set.

First we recall definition of  a directional solution set.
\begin{defn}[Directional Solution Set]\cite[Definition 4.5]{BYZ}\label{Defn4.1}
	Let $\bar x, u\in \mathbb{R}^n$. The set of optimal solutions for problem $(P_x)$ in direction $u$ is defined by
	$$S(\bar x;u)=\{y\in S(\bar x)| \exists x^k\xrightarrow{u}\bar x,
	 y^k\rightarrow y, \ \mbox{s.t.} \ y^k\in S(x^k)\}.$$  
\end{defn}
It is easy to see that for a nonzero direction $u$, a directional solution set is a subset of the solution set. If the directional solution set is nonempty, then we may try to use it to replace the solution set and provide tighter bounds to
 the upper/lower Dini directional derivative of the value function
 in (\ref{minmax}).
{Generally, it is possible that $S(\bar x;u)=\emptyset$ even when the solution set is nonempty and so the generalization using the directional solution set is only useful 
  when the directional solution is nonempty. To answer this question, we first present the following lemma.
 \begin{lemma}\label{lemma4.1}Let $\bar x, u\in \mathbb{R}^n$. Suppose that  $V(x)$ is continuous at $\bar x$ in direction $u$. Let $t_k \downarrow 0, u^k\rightarrow u$ and $y^k \in S(\bar x+t_k u^k)$. If $y_k\rightarrow  y$, then $ y\in S(\bar x;u)$.
 \end{lemma}
 \beginproof Since $y^k \in S(\bar x+t_k u^k)$, we have 
\begin{eqnarray*}
	&& f(\bar x, y)=\lim_{k\rightarrow \infty}f(\bar x+t_ku^k,y^k)=\lim_{k\rightarrow \infty}{ V}(\bar x+t_ku^k)= { V}(\bar x),\\
	&& P(\bar x, y)= \lim_{k\rightarrow \infty} P(\bar x+t_ku^k,y^k) \in C,
\end{eqnarray*}  which implies that $ y\in S(\bar x)$. By Definition \ref{Defn4.1}, we have
	$ y\in S(\bar x;u).
	$
 \endproof
 We now give a sufficient condition for the nonemptyness of the directional solution set.
 First we present a condition which guarantees the directional lower semicontinuity of the value function. 
\begin{defn}[Directional Restricted Inf-compactness]\cite[Definition 4.1]{BY}
	We say that the restricted inf-compactness holds at $\bar x$ in direction $u$ if $V(\bar x)$ is finite and there exist a compact set  $\Omega_u\subseteq \mathbb{R}^n$, and positive numbers $\varepsilon>0, \delta>0$ such that for all $ x\in \bar x+{\cal V}_{\varepsilon,\delta}(u)$ with $V(x)<V(\bar x)+\varepsilon$, one always has $S(x)\cap\Omega_u\neq\emptyset$.
\end{defn}
Obviously, if the restricted inf-compactness holds at $\bar x$ in {direction u=0}, then the {classical} restricted inf-compactness defined in \cite[Hypothesis 6.5.1]{Clarke} and \cite[Definition 3.8]{GLYZ} holds. The reader is referred to some easier to verify sufficient conditions for the directional restricted inf-compactness in \cite{BY}.

We now show that under the directional restricted inf-compactness condition and directional continuity of the value function, the directional solution set is nonempty.
\begin{prop}\label{directsolution} Suppose the restricted inf-compactness condition holds at $\bar x$ in direction $u$. If $V(x)$ is continuous at $\bar x$ in direction $u$. Then $S(\bar x; u)$ is nonempty.
\end{prop}
\beginproof
Let $t_k\downarrow0, u^k\rightarrow u$.
 Since $V(x)$ is continuous at $\bar x$ in direction $u$,  for any $\varepsilon>0$  $V(\bar x+t_ku^k)< V(\bar x)+\varepsilon$ for sufficiently large $k$. Since the restricted inf-compactness holds at $\bar x$ in direction $u$ with the compact set $\Omega_u$, there exists a sequence $y^k\in S(\bar x+t_ku^k)\cap\Omega_u$ for $k$ large enough. Without loss of generality, assume $y:=\lim_{k\rightarrow \infty}y^k$. Then by Lemma \ref{lemma4.1}, 
	$ y\in S(\bar x;u)\cap\Omega_u.
	$
	\endproof
	
	Now we give a sufficient condition for the directional continuity of the value function which is a directional version of the result in Guo et al. \cite[Proposition 3.1]{GYZ2021}. 



{
\begin{prop}
Assume that the restricted inf-compactness holds at $\bar x$ in direction $u$ and there exists $\bar y\in S(\bar x)$ such that 
\begin{equation}\displaystyle \lim_{x\xrightarrow{u} \bar{x}} {\rm dist}(\bar{y},\mathcal F(x)) =0. \label{eqn(27)} \end{equation}
Then $V(x)$ is continuous at $\bar x$ in direction $u$.
\end{prop}}
\beginproof Under the restricted inf-compactness at $\bar x$ in direction $u$, it is easy to verify that $V(x)$ is l.s.c. at $\bar x$ in direction $u$ which can be proved using  a similar argument as in the proof of the non-directional case in Guo et al. \cite[Theorem 3.9]{GLYZ}.  We now prove the upper semicontinuity of $V(x)$ at $\bar x$ in direction $u$. Let $\{x^k\}$ be a sequence such that $x^k\xrightarrow{u} \bar{x}$ and 
$$\lim_{k\rightarrow \infty} V(x^k)= \displaystyle \limsup_{x\xrightarrow{u} \bar{x}}
		V(x).$$ Then 
		$\displaystyle \lim_{k\to\infty} {\rm dist}(\bar{y},\mathcal F(x^k))=0.$ For each $k$, take $y^k\in \mathcal F(x^k)$ such that
	$\|y^k-\bar{y}\| < {\rm dist}(\bar{y},\mathcal F(x^k)) + 1/k$. It follows that $ y^k \to \bar{y}$ and so
	\[
	V(\bar{x}) = f(\bar{x}, \bar{y}) = \lim_{k\to\infty}f(x^k, y^k) \ge \lim_{k\to\infty}V(x^k)
	= \limsup_{x\xrightarrow{u}  \bar{x}}V(x).
	\]
	This shows that $V(x)$ is u.s.c. at $\bar{x}$ in direction $u$. 
	The proof is complete.
\endproof

It is obvious that the RS at $(\bar x,\bar y)$ in direction $u$ for some $\bar y\in S(\bar x)$ implies condition (\ref{eqn(27)}). Hence we have the following corollary.
\begin{cor}
	Assume that the restricted inf-compactness holds at $\bar x$ in direction $u$ and the feasible map ${\cal F}$ satisfies RS at $(\bar x,\bar y)$ in direction $u$ for some $\bar y\in S(\bar x)$, then $V(x) $ is continuous at $\bar x$ in direction $u$ and $S(\bar x;u)\not =\emptyset$.
\end{cor}
%

\begin{defn}[Directional Inner Semicontinuity]\cite[Definition 4.3]{BY} Given $\bar y\in S(\bar x)$, we say that  the optimal solution map $S(x)$ is inner semicontinuous at $(\bar x,\bar y)$ in direction $u$, if for any sequence $x^k \xrightarrow{u}  \bar{x}$, there exists a sequence $y^k\in S(x^k)$ converging to $\bar 
	y$. 
\end{defn}
By definition, if $\exists \bar y\in S(\bar x)$ such that $S(x)$ is inner semicontinuous  at $(\bar x,\bar y)$ in direction $u$, then the restricted inf-compactness holds at $\bar x$ in direction $u$. Note that if $S(x)$ is inner semicontinuous at $(\bar x,\bar y)$ in {direction u=0}, then $S(x)$ is inner semicontinuous at $(\bar x,\bar y)$ in the sense of \cite[Definition 1.63]{Aub2}. When the solution set $S(x)$ is a single-valued map, the directional inner semicontinuity is reduced to the directional continuity of the mapping $S(x)$. 
As commented in \cite{BY},  the directional inner semicontinuity of $S(x)$ at $(\bar x,\bar y)$ in direction $u$ implies the restricted inf-compactness condition and hence the lower semicontinuity of $V(x)$ at $\bar x$ in direction $u$. Moreover if $S(x)$ is inner semicontinuous at $(\bar x,\bar y)$ in direction $u$, then
$${\rm dist} (\bar y, {\cal F}(x_k)) \leq {\rm dist} (\bar y, S(x_k)) \leq \|\bar y-y^k\|$$ 
and hence condition (\ref{eqn(27)}) must hold and so the value function $V(x)$ must be continuous. Moreover from definition, $S(x)$ is inner semicontinuous at $(\bar x,\bar y)$ in direction $u$ implies that $\bar y\in S(\bar x;u)$.
Hence we have the following corollary.
\begin{cor}\label{cor4.2} Suppose that $S(x)$ is inner semicontinuous at $(\bar x,\bar y)$ in direction $u$, then  $V(x)$ is continuous at $\bar x$ in direction $u$ and $\bar y\in S(\bar x;u)$.
\end{cor}

%

 We now extend \cite[Theorem 4.26]{BS} to allow a nonconvex set $C$.  It is obvious that the regularity condition (\ref{wRCQ4}) is in general weaker than the Robinson's CQ (\ref{RCQ}).

\begin{thm}\label{nondir} Assume that the restricted inf-compactness holds at $\bar x$ in direction $u$ with compact set $\Omega_u$  and RS holds for the problem $(P_{\bar x})$ 
at each $y\in {S(\bar x)\cap\Omega_u}$ in direction $u$. 
Suppose  for  each $y\in {S(\bar x)\cap\Omega_u}$ the conic linear system in variable $v$:
\begin{equation}\label{feasibleso}
	\nabla P(\bar x,y)(\pm u, v)\in \widehat T_C(P(\bar x,y))
\end{equation}
is feasible and the  regularity condition
\begin{equation}\label{wRCQ4}
	0=\nabla_y P(\bar x,y)^T\lambda,\ \lambda\in N_C^c (P(\bar x,y))\cap \{\nabla P(\bar x,y)(\pm u, v)\}^\perp \ \implies \lambda=0
\end{equation} holds for any $v$ satisfying system (\ref{feasibleso}).
Then \begin{align*}
  {\min_{ y\in S(\bar x;u)\cap\Omega_u}}
 {\min_{\lambda\in  \Lambda^c(\bar x, y)}}\nabla_x L(\bar x,y, \lambda)u\leq V_-'(\bar x;u) \leq V_+'(\bar x;u)\leq{\min_{ y\in S(\bar x;u)\cap\Omega_u}}
	 {\max_{\lambda\in  \Lambda^c(\bar x, y)}}\nabla_x L(\bar x,y, \lambda)u.
\end{align*} 
Furthermore if  $\Lambda^c(\bar x, y)=\{\bar \lambda(y)\}$ is a singleton for every 
$y\in S(\bar x;u)\cap\Omega_u$, 
then $V(x)$ is Hadamard directionally differentiable and 
$$V'(\bar x; u)= \min_{ y\in S(\bar x;u)\cap\Omega_u}
\nabla_x L(\bar x, y, \bar \lambda(y))u.$$

\end{thm}
\beginproof
Recall that for any closed convex cone $\cal K$ and any $d\in \cal K$, 
$N_{\cal K}(d)={\cal K}^o\cap\{d\}^\perp$; see e.g.  \cite[Corollary 23.5.4]{CA}.
Since the regular tangent cone $\widehat{T}_C(P(\bar x,y))$ is  closed and convex, 
\begin{align*}
	N_{\widehat{T}_C(P(\bar x,y))}(\nabla P(\bar x,y)(\pm u,v))&= \{\widehat{T}_C(P(\bar x,y))\}^\circ\cap\{\nabla P(\bar x,y)(\pm u,v)\}^\perp\\ &=  N^c_C(P(\bar x,y))\cap\{\nabla P(\bar x,y)(\pm u,v)\}^\perp,
\end{align*}
where the second equality follows from tangent-normal polarity in Proposition \ref{polarity}. 
{It follows that (\ref{wRCQ4}) is equivalent to that Robinson's CQ for the conic linear system (\ref{feasibleso}) holds at each solution $v$, i.e.,
\begin{equation}\label{wRCQ3}
	0=\nabla_y P(\bar x,y)^T\lambda,\ \lambda\in N_{\widehat{T}_C(P(\bar x,y))}(\nabla P(\bar x,y)(\pm u, v)) \ \implies \lambda=0.
\end{equation}
By the discussion before Theorem \ref{dt}, feasibility of the conic linear system (\ref{feasibleso}) and condition (\ref{wRCQ3}) at each feasible solution $v$ is equivalent to the following condition
\begin{equation}\label{wRCQ}
	0\in{\rm int} \left \{{\pm }\nabla_x P(\bar x,y)u+ \nabla_yP(\bar x,y)\mathbb R^m+\widehat{T}_C(P(\bar x,y)) \right \}.
\end{equation}}

Now we prove the estimates.
By \cite[Theorem 4.1]{BY}, since the restricted inf-compactness in direction $u$ and RS holds at $(\bar x,y)$ for each $y\in S(\bar x)$, $V(x)$ is directionally Lipschitz continuous at $\bar x$ in direction $u$, hence by Proposition \ref{directsolution} $S(\bar x;u)\not =\emptyset$. Moreover 
$$-\infty < V'_-(\bar x;u) \leq V'_+(\bar x;u) <\infty.$$

{\bf (I) Proof of the lower bound}.
Let $t_k\downarrow0, u^k\rightarrow u$ be the sequences satisfying
\begin{equation*}
	V'_-(\bar x;u)=\lim_{k\rightarrow \infty}\frac{V(\bar x+ t_k  u^k)-V(\bar x)}{ t_k}.
\end{equation*} By the continuity of $V$ at $\bar x$ in direction $u$, for any $\varepsilon>0$, $V(\bar x+t_k u^k) < V(\bar x)+\varepsilon$ for sufficiently large $k$.  By the restricted inf-compactness condition at $\bar x$ in direction $u$,
there exists a sequence $y^k\in S(\bar x+t_ku^k)\cap\Omega_u$ for $k$ large enough. Without loss of generality, assume $\tilde y:=\lim_{k\rightarrow \infty}y^k$. Then by Lemma \ref{lemma4.1}, we have $\tilde y\in S(\bar x;u)\cap\Omega_u$.
	Since {condition (\ref{wRCQ})} holds at $(\bar x,\tilde y)$, there exists $\tilde v\in \mathbb{R}^m$ such that 
	\begin{equation} -\nabla P(\bar x,\tilde y)(u,\tilde v) \in \widehat T_C(P(\bar x,\tilde y)). \label{choice} \end{equation}
	Now consider any vector $\tilde v$ satisfying (\ref{choice}). By definition of the regular tangent cone $\widehat T_C(P(\bar x,\tilde y))$ in Definition \ref{tangentnormal}, since $P(\bar x+t_ku^k,y^k) \rightarrow P(\bar x,\tilde y)$ as $k\rightarrow \infty$,
	we can choose $  s^k\in C  $ for each $k$ such that 
	\begin{equation}\label{s^k}
		\lim_{k}\frac{s^k-P(\bar x+t_ku^k,y^k)}{t_k}=-\nabla P(\bar x, \tilde y)(u,\tilde v).
	\end{equation}
	For any $\varepsilon>0$, since $P(x,y)$ is continuously differentiable, we always have \begin{equation} \label{epsilon}
	\|\nabla P(\bar x+ t_k u, y^k)-\nabla P(\bar x,\tilde y)\|\leq\varepsilon
	\end{equation} for sufficiently large $k$. Then by the RS  at $(\bar x,\tilde y)$ in direction $u$,  there exists $\kappa>0$ such that
	\begin{align*}
		\mbox{dist}(y^k-t_k\tilde v,\mathcal F(\bar  x))&\leq\kappa\mbox{dist}(P(\bar x,y^k-t_k\tilde v),C)\\
		&=\kappa\mbox{dist}(P(\bar x+t_ku, y^k)- t_k\nabla P(\bar x+ t_ku, y^k)(u,\tilde v)+o(t_k),C)\\
		&\leq\kappa\|P(\bar x+t_ku,y^k)-t_k\nabla P(\bar x+t_ku, y^k)(u,\tilde v)+o(t_k)-s^k\|\\
		&\leq\kappa\|P(\bar x+t_ku,y^k)-t_k\nabla P(\bar x,\tilde y)(u,\tilde v)+o(t_k)-s^k\|\\
		&\quad +\kappa t_k\|-\nabla P(\bar x+t_ku, y^k)+\nabla P(\bar x,\tilde y)\|\|(u,\tilde v)\|\\
		&=O(t_k\varepsilon),
	\end{align*}
	where the last equality follows from (\ref{s^k}) and (\ref{epsilon}). This means for each $k$ sufficiently large there exists $\hat{y}^k\in \mathcal F(\bar x)$ such that $\|\hat y^k-(y^k-t_k\tilde v)\|=O(t_k\varepsilon)$ and so $y^k=\hat y^k+t_k \tilde v+O(t_k\varepsilon) {\bf e}$,  where ${\bf e}$ is a unit vector.
	Consequently, we have $\lim_k\hat y^k=\tilde y$ and
	\begin{eqnarray*}
		\nonumber
		V'_-(\bar x;u)=\lim_{k\rightarrow \infty}\frac{V(\bar x+ t_k u^k)-V(\bar x)}{t_k}
		&& {=} \lim_{k\rightarrow\infty}\frac{f(\bar x+t_k u^k, y^k)-f(\bar x,\tilde y)}{t_k}\\
		&& \geq\lim_{k\rightarrow\infty}\frac{f(\bar x+t_k u^k, y^k)-f(\bar x, \hat y^k)}{t_k}\\
		\nonumber
		&&=\lim_{k\rightarrow\infty}\frac{f(\bar x+t_ku^k, \hat y^k+t_k \tilde v+{O(t_k\varepsilon) {\bf e}})-f(\bar x,\hat y^k)}{t_k}\\
		&&=\nabla f(\bar x,\tilde y)(u,\tilde v)+O(\varepsilon).
	\end{eqnarray*}
	Since this is true for any $\varepsilon>0$, we obtain 
	\begin{equation}\label{eqn(30)}
		V_-'(\bar x;u)\geq\nabla f(\bar x,\tilde y)(u, \tilde v),\qquad  \mbox{for any } \tilde v \mbox{ satisfying  (\ref{choice})}.
	\end{equation}
	Thus 
	\begin{eqnarray*}
		V_-'(\bar x;u) 
		&\geq  &\sup_{\tilde v} \left \{\nabla f(\bar x,\tilde y)(u,\tilde v): -\nabla P(\bar x,\tilde y)(u, \tilde v)\in \widehat{T}_C(P(\bar x,\tilde y))\right \}\nonumber\\
		& =&- \inf_{\tilde v} \left \{\nabla f(\bar x,\tilde y)(-u,-\tilde v): \nabla P(\bar x,\tilde y)(-u, -\tilde v)\in \widehat{T}_C(P(\bar x,\tilde y))\right \}.
	\end{eqnarray*}
%
Since the regular tangent cone is closed and convex, the following problem is  a conic linear program in the form of 
(\ref{primalP}).
\begin{align*}
(P_u)~~~~~~~	&\min_{\tilde v}   \nabla f(\bar x,\tilde y)(-u,-\tilde v) \ \nonumber\\
	& s.t. ~ \nabla P(\bar x,\tilde y)(-u, -\tilde v)\in \widehat{T}_C(P(\bar x,\tilde y)).
\end{align*}
 By virtue of condition (\ref{wRCQ}), problem $(P_u)$  is feasible. 
 By (\ref{dualD}), 
the Lagrange dual program of $(P_u)$ can be written in the form
\begin{align*}
(D_u)~~~~~	&
	&\max_{\lambda\in (\widehat{T}_C(P(\bar x,\tilde y))^\circ}\ &\nabla_x L(\bar x,\tilde y, \lambda)(-u)\notag\\
	& &s.t.\quad\quad~ &\nabla_y f(\bar x,\tilde y)+\nabla_y P(\bar x,\tilde y)^T\lambda=0.
\end{align*}

By the regularity condition (\ref{wRCQ}), program {\rm (P$_u$)} is feasible.
Therefore val$(P_u)<+\infty$. Since $V'_-(\bar x;u)<\infty$ and (\ref{eqn(30)}) holds, we have val$(P_u)>-\infty$. 
By Theorem \ref{dt} the strong duality of {\rm (P$_u$)} then follows with a nonempty solution set of the dual problem.
  By Proposition \ref{polarity}, $(\widehat{T}_C(P(\bar x,\tilde y))^\circ=N^c_C(P(\bar x,\tilde y))$ and hence
\begin{eqnarray*}
		V_-'(\bar x;u) 
		& \geq &- \inf_{\tilde v} \left \{\nabla f(\bar x,\tilde y)(-u,-\tilde v): \nabla P(\bar x,\tilde y)(-u, -\tilde v)\in \widehat{T}_C(P(\bar x,\tilde y))\right \}\\
	&=& -\max_{\lambda\in  \Lambda^c(\bar x, \tilde y)}\nabla_x L(\bar x, \tilde y,\lambda)(-u)=\min_{\lambda\in  \Lambda^c(\bar x, \tilde y)}\nabla_x L(\bar x, \tilde y,\lambda)u.
\end{eqnarray*}
{\bf (II) Proof of the upper bound}.
Let $t_k\downarrow0, u^k\rightarrow u$ be the sequences satisfying
\begin{equation*}
	V'_+(\bar x;u)=\lim_{k\rightarrow \infty}\frac{V(\bar x+ t_k  u^k)-V(\bar x)}{ t_k}.
\end{equation*}
Take any $\tilde y\in {S(\bar x;u)\cap\Omega_u}$. 
For any direction $v$ satisfying that $\nabla P(\bar x,\tilde y)(u,v)\in \widehat T_C(P(\bar x,\tilde y))$,  by definition of the regular tangent cone in Definition \ref{tangentnormal}, {for the sequences 
$t_k\downarrow0$, there exists $s^k\in C$ such that} 
$$
\lim_k\frac{s^k-P(\bar x,\tilde y)}{t_k}=\nabla P(\bar x,\tilde y)(u,v).
$$
Since the RS holds at $(\bar x, \tilde y)$ in direction $u$, we have
\begin{align*}
\mbox{dist}(\tilde y+t_kv,\mathcal F(\bar x+t_ku^k))&\leq\kappa\mbox{dist}(P(\bar x+t_ku^k,\tilde y+t_kv),C)\\
&\leq\kappa\|P(\bar x+t_ku^k,\tilde y+t_kv)-s^k\|\\
&=\kappa\|P(\bar x+t_ku^k,\tilde y+t_kv)-(P(\bar x,\tilde y)+t_k\nabla P(\bar x,\tilde y)(u,v)+o(t_k))\|\\
&=o(t_k).
\end{align*}
This means that there exists a sequence $\hat y^k\in \mathcal F(\bar x+t_ku^k)$ such that $\hat y^k=\tilde y+t_kv+o(t_k)$, i.e., $(\hat y^k-\tilde y)/t_k\rightarrow v$.

 Then we obtain
\[
\lim_k\frac{V(\bar x+t_ku^k)-V(\bar x)}{t_k}
\leq\lim_k\frac{f(\bar x+ t_ku^k,\hat y^k)-f(\bar x,\tilde y)}{t_k}=\nabla f(\bar x,\tilde y)(u,v).
\]
Since $\tilde y$ and $v$ are arbitrary,  we have 
\begin{eqnarray*}
	V'_+(\bar x;u) &\leq&  \inf_v \{ \nabla f(\bar x,\tilde y)(u,v):\nabla P(\bar x,\tilde y)(u,v)\in \widehat T_C(P(\bar x,\tilde y)) \}.
	\end{eqnarray*}
Similarly as the proof of the lower bound, by using the strong duality,
we have \begin{eqnarray*}
V'_+(\bar x;u)	&\leq & \max_{\lambda\in  \Lambda^c(\bar x, \tilde y)}\nabla_x L(\bar x, \tilde y,\lambda)u.
\end{eqnarray*}
It follows that
\begin{eqnarray*}
	V'_+(\bar x;u)\leq{\min_{ y\in S(\bar x;u)\cap\Omega_u}}\max_{\lambda\in  \Lambda^c(\bar x, y)}\nabla_x L(\bar x, y,\lambda)u.
\end{eqnarray*}
\endproof

\begin{remark}\label{wRCQl} 
When set $C$ is a  closed convex set, we have
\begin{eqnarray*}
 N_C(P(\bar x,y))\cap\{\nabla P(\bar x,y)(\pm u,-v)\} ^\perp
&=& N_C(P(\bar x,y);\nabla P(\bar x,y)(\pm u,-v)).
\end{eqnarray*} 
So if $C$ is convex, then the regularity condition (\ref{wRCQ4}) is equivalent to the FOSCMS  in direction $\nabla P(\bar x,y)(\pm u,-v)$:
\begin{equation*}
	0\in\nabla_y P(\bar x,y)^T\lambda,\ \lambda\in N_C(P(\bar x,y);\nabla P(\bar x,y)(\pm u,-v)) \ \implies \lambda=0.
\end{equation*}
{If the set $C$ is convex, then the regular tangent cone coincides with the tangent cone and hence the existence of $v$ satisfying system (\ref{feasibleso}) is reduced to  the existence of $v$ satisfying system
$$ \nabla P(\bar x,y)(\pm u, v)\in  T_C(P(\bar x,y)).$$
The feasibility of the above system can be guaranteed by RS in direction $\pm u$ as demonstrated in Proposition \ref{LinearizC}.}
\end{remark}

By Proposition \ref{prop3.2}, if $P(x,y)=a(x)+By+c$ with $a:\mathbb R^n\rightarrow\mathbb R^p$ being continuous, $B\in\mathbb R^{p\times m}, c\in \mathbb{R}$, $C$ is the union of finitely many convex polyhedral sets and {${\cal F}(x)$ is nonempty in a neighborhood of $\bar x$ in direction $u$}, then RS in direction $u$ holds at each $(\bar x,\bar y)\in {\rm gph} {\cal F}$. {Moreover when $C$ is the union of finitely many convex polyhedral sets, $\widehat{T}_C(P(\bar x, y))$ is a convex polyhedral set. Moreover Remark \ref{wRCQl}, if $C$ is a convex polyhedral set, then RS implies the feasibility of the linearized system. Hence by the discussion before Theorem \ref{dt}, the strong duality holds for the conic linear program $(P_u)$ without imposing the regularity condition (\ref{wRCQ4}).}
 Hence from the proof of Theorem \ref{nondir} and Remark \ref{wRCQl}, we have the following result for this special case.
\begin{thm} \label{cor4.3} Assume that the restricted inf-compactness holds at $\bar x$ in direction $u$ with compact set $\Omega_u$.  Suppose that   $P(x,y)=a(x)+By+c$ with $a:\mathbb R^n\rightarrow\mathbb R^p$ being continuously differentiable, $B\in\mathbb R^{p\times m}, c\in \mathbb{R}$, $C$ is the union of finitely many convex polyhedral sets and {${\cal F}(x)$ is nonempty in a neighborhood of $\bar x$ in direction $u$}.
Suppose that either $C$ is a convex polyhedral set or for each $y\in S(\bar x)$, 
\begin{equation*} \exists v \mbox{ such that } \qquad \nabla P(\bar x, \bar y)(\pm u,v)\in \widehat{T}_C(P(\bar x, \bar y)).
\end{equation*}
Then 
$$ {\min_{ y\in S(\bar x;u)\cap\Omega_u}}
 \min_{\lambda\in  \Lambda^c(\bar x, y)}\nabla_x L(\bar x,y, \lambda)u\leq V_-'(\bar x;u) \leq V_+'(\bar x;u)
	\leq{\min_{ y\in S(\bar x;u)\cap\Omega_u}}
	 \max_{\lambda\in  \Lambda^c(\bar x, y)}\nabla_x L(\bar x,y, \lambda)u.$$
	 Furthermore  if  $\Lambda^c(\bar x, y)=\{\bar \lambda(y)\}$ is a singleton for every 
$y\in S(\bar x;u)\cap\Omega_u$,
then $V(x)$ is Hadamard directionally differentiable  at $\bar x$ in direction $u$ and 
$$V'(\bar x; u)= {\min_{ y\in S(\bar x;u)\cap\Omega_u}}
\nabla_x L(\bar x, y, \bar \lambda(y))u.$$
\end{thm}

When the solution map has  directional inner semicontinuity, the estimates of the directional derivative is simpler using only one solution as shown in the corollary below.
\begin{cor} \label{isc} Suppose that the solution set $S(\bar x)$ is inner semicontinuous at $(\bar x,\bar y) \in {\rm gph} S$ in direction $u$ and RS holds for the problem $(P_{\bar x})$ 
at $\bar y$ in direction $u$. 
Suppose   the conic linear system in variable $v$:
	$\nabla P(\bar x,\bar y)(\pm u, v)\in \widehat T_C(P(\bar x,\bar y))$
is feasible and for any feasible $v$ satisfying the above system, the regularity condition
	$$0=\nabla_y P(\bar x,\bar y)^T\lambda,\ \lambda\in N_C^c (P(\bar x,\bar y))\cap \{\nabla P(\bar x,\bar y)(\pm u, v)\}^\perp \ \implies \lambda=0$$
 holds.
Then \begin{align*}
  {}
 {\min_{\lambda\in  \Lambda^c(\bar x, \bar y)}}\nabla_x L(\bar x, \bar y, \lambda)u\leq V_-'(\bar x;u) \leq V_+'(\bar x;u)\leq{}
	 {\max_{\lambda\in  \Lambda^c(\bar x, \bar y)}}\nabla_x L(\bar x,  \bar y, \lambda )u.
\end{align*} 
Furthermore if  $\Lambda^c(\bar x, \bar y)=\{\bar \lambda\}$ is a singleton, 
then $V(x)$ is Hadamard directionally differentiable and 
$$V'(\bar x; u)= 
\nabla_x L(\bar x,\bar y, \bar \lambda)u.$$
\end{cor}
\beginproof
Since $S(x)$ is inner semicontinuous at $(\bar x,\bar y)\in{\rm gph}S$ in direction $u$, for any sequences $t_k\downarrow0, u^k\rightarrow u$ there exists a sequence $y^k\in S(\bar x+t_ku^k)$ converging to $\bar y$. Then the restricted inf-compactness holds at $\bar x$ in direction $u$ with $\Omega_u:=\bar y+\epsilon\mathbb B$ for any $\epsilon>0$. By Corollary \ref{cor4.2}, $\bar y\in S(\bar x;u)$. And we can always choose $\tilde y=\bar y$ in the proof of Theorem \ref{nondir}. Hence the result follows from Theorem 
\ref{nondir}. 
\endproof

Similarly by Theorem \ref{cor4.3} and Corollary  \ref{isc} we have the following results for the special case where $P(x,y)$ is affine in variable $y$ and $C$ is the union of finitely many convex  polyhedral sets.
\begin{cor}\label{cor-nodir} Suppose that the solution set $S(x)$ is inner semicontinuous at $(\bar x,\bar y) \in {\rm gph} S$ in direction $u$. Then $\bar y\in S(\bar x;u)$. If  $P(x,y)=a(x)+By+c$ with $a:\mathbb R^n\rightarrow\mathbb R^p$ being continuously differentiable, $B\in\mathbb R^{p\times m}, c\in \mathbb{R}$, $C$ is the union of finitely many convex  polyhedral sets and  {${\cal F}(x)$ is nonempty in a neighborhood of $\bar x$ in direction $u$}. If either $C$ is a convex polyhedral set or
$$ \exists v \mbox{ such that } \qquad \nabla P(\bar x, \bar y)(\pm u,v)\in \widehat{T}_C(P(\bar x, \bar y)),$$ then 
$$ \min_{\lambda\in  \Lambda^c(\bar x, \bar y)} \nabla_x L(\bar x,\bar y,\lambda)u\leq V_-'(\bar x;u) \leq V_+'(\bar x;u)\leq 
	\max_{\lambda\in  \Lambda^c(\bar x, \bar y)}\nabla_x L(\bar x,\bar y, \lambda)u.$$
	 Furthermore if  $\Lambda^c(\bar x, \bar  y)=\{ \bar \lambda \}$ is a singleton, 
then $V(x)$ is Hadamard directionally differentiable  at $\bar x$ in direction $u$ and 
$V'(\bar x; u)= 
\nabla_x L(\bar x,\bar   y, \bar \lambda)u.$
\end{cor}

By Corollary \ref{cor3.1}, we know that when $\nabla_xP(\bar x, y)u\neq 0$,  FOSCMS in  directions
$\nabla P(\bar x, y)(\alpha u, v)$ with $\alpha \geq 0$ can be used to verify directional RS in direction $u$. Since the regularity condition (\ref{wRCQ4}) is  FOSCMS in  directions $\nabla P(\bar x, y)(\pm u, v)$, we can strengthen the regularity condition (\ref{wRCQ4}) to cover First-order sufficient conditions for  RS in direction $u$ required in Corollary \ref{cor3.1}.
By Corollary \ref{cor3.1}, Theorem \ref{nondir} and Remark \ref{wRCQl}, we have the following corollary immediately. 
{\begin{cor} Assume that $C$ is closed and convex and the restricted inf-compactness holds at $\bar x$ in direction $u$ with compact set $\Omega_u$. Suppose for each $ y\in S(\bar x)\cap \Omega_u$,
$\nabla_xP(\bar x, y)u\neq 0$ and for every $\alpha  \in \mathbb{R}$
and every $\tau_k\downarrow 0$,  
there exists $v$ such that 
		\begin{equation*}
			\liminf_{k\rightarrow \infty }
			\frac{{\rm dist}(P(\bar x, y)+\tau_k\nabla P(\bar x, y)(\alpha u, v),C)}{\tau_k}
			=0.
		\end{equation*}
		If for every $ y\in S(\bar x)\cap \Omega_u$, $(0,0)\neq(\alpha,v)\in \mathbb{R}\times \mathbb{R}^m $ with $\nabla P(\bar x, y)(\alpha u, v) \in T_C(P(\bar x, y))$ one has 
		\begin{equation*}
			\nabla_y P(\bar x, y)^T\lambda=0,\ \lambda\in N_C(P(\bar x, y))\cap 
			\{\nabla P(\bar x, y)(\alpha u, v )\}^\perp
			\implies\lambda=0. 
		\end{equation*} 
		Then \begin{align*}
  {\min_{ y\in S(\bar x;u)\cap\Omega_u}}
 {\min_{\lambda\in  \Lambda(\bar x, y)}}\nabla_x L(\bar x,y, \lambda)u\leq V_-'(\bar x;u) \leq V_+'(\bar x;u)\leq{\min_{ y\in S(\bar x;u)\cap\Omega_u}}
	 {\max_{\lambda\in  \Lambda(\bar x, y)}}\nabla_x L(\bar x,y, \lambda)u.
\end{align*} 
Furthermore if  $\Lambda(\bar x, y)=\{\bar \lambda(y)\}$ is a singleton for every 
$y\in S(\bar x;u)\cap\Omega_u$, 
then $V(x)$ is Hadamard directionally differentiable and 
$$V'(\bar x; u)= \min_{ y\in S(\bar x;u)\cap\Omega_u}
\nabla_x L(\bar x, y, \bar \lambda(y))u.$$
\end{cor}
\beginproof By Corollary \ref{cor3.1}, ${\cal F}$ satisfies RS at $(\bar x,\bar y)$ in direction $\pm u$. The rest of the proof follows from Theorem \ref{nondir} and Remark \ref{wRCQl}.
\endproof}
Similarly by  Corollaries \ref{cor3.1} and \ref{isc},  and Remark \ref{wRCQl}, we have the following result.
\begin{cor} \label{isc2} Assume that $C$ is closed and convex and  the solution set $S(\bar x)$ is inner semicontinuous at $(\bar x,\bar y) \in {\rm gph} S$ in direction $u$.
Suppose that $\nabla_xP(\bar x, \bar y)u\neq 0$ and for every {$\alpha  \in \mathbb{R}$}
and every $\tau_k\downarrow 0$,  
there exists $v$ such that 
		\begin{equation*}
			\liminf_{k\rightarrow \infty }
			\frac{{\rm dist}(P(\bar x, \bar y)+\tau_k\nabla P(\bar x, \bar y)(\alpha u, v),C)}{\tau_k}
			=0.
		\end{equation*} If for any $(0,0)\neq(\alpha,v)\in \mathbb{R}\times \mathbb{R}^m $ with $\nabla P(\bar x, \bar y)(\alpha u, v) \in T_C(P(\bar x, \bar y))$ one has 
		\begin{equation*}
			\nabla_y P(\bar x, \bar y)^T\lambda=0,\ \lambda\in N_C(P(\bar x, \bar y))\cap 
			\{\nabla P(\bar x, \bar y)(\alpha u, v )\}^\perp
			\implies\lambda=0. 
		\end{equation*}  
Then \begin{align*}
 {\min_{\lambda\in  \Lambda^c(\bar x, \bar y)}}\nabla_x L(\bar x,\bar    y, \lambda )u\leq V_-'(\bar x;u) \leq V_+'(\bar x;u)\leq{}
	 {\max_{\lambda\in  \Lambda^c(\bar x, \bar y)}}\nabla_x L(\bar x,  \bar y, \lambda )u.
\end{align*} 
\end{cor}


If Robinson's CQ holds for the problem ({\rm P$_{\bar x}$}) at each $y\in {S(\bar x;u)\cap\Omega_u}$,  then  the condition (\ref{Robinsonnew}) holds which implies that  (\ref{wRCQ}) holds automatically and each multiplier set $\Lambda^c(\bar x,  y)$ is compact. Moreover by Proposition \ref{RCQRS}, the RS holds automatically.  Hence the following corollary follows immediately.
\begin{cor} \label{cor4.4} Suppose that the restricted inf-compactness holds at $\bar x$ in direction $u$ with compact set $\Omega_u$ and  Robinson's CQ holds for problem $(P_{\bar x})$ at each $ y\in S(\bar x;u)$.  Then 
$$ {\min_{ y\in S(\bar x;u)\cap\Omega_u}}
 {\min_{\lambda\in  \Lambda^c(\bar x, y)}}\nabla_x L(\bar x,y, \lambda)u\leq V_-'(\bar x;u) \leq V_+'(\bar x;u)
	\leq{\min_{ y\in S(\bar x;u)\cap\Omega_u}}
	 {\max_{\lambda\in  \Lambda^c(\bar x, y)}}\nabla_x L(\bar x,y, \lambda)u.$$
	 Furthermore if  $\Lambda^c(\bar x, y)=\{\bar \lambda(y)\}$ is a singleton for every 
$y\in S(\bar x;u)\cap\Omega_u$, 
then $V(x)$ is Hadamard directionally differentiable  at $\bar x$ in direction $u$ and 
$$V'(\bar x; u)= {\min_{ y\in S(\bar x;u)\cap\Omega_u}}
\nabla_x L(\bar x,y, \bar \lambda(y))u.$$
\end{cor}
Even in the case when $C$ is convex and Robinson's CQ holds,  our result in Corollary \ref{cor4.4} provides  sharper results than 
\cite[Theorem 4.26]{BS}. In particular to conclude Hadamard directional differentiability, we do not need to assume that the set of multipliers is a singleton for every $y\in S({\bar x})$. Instead we only need to assume that the set of multipliers is a singleton for every  $y$ in the directional solution set $S(\bar x;u)$. Our assumptions are weaker and the bounds are tighter  since  the directional solution set $S(\bar x; u)$ is in general smaller than the solution set $S(\bar x)$.

{
Finally, we state the special case of parametric nonlinear program where $P(x,y):=(h(x,y),g(x,y))$, $C:=\{0\}_{p_1}\times\mathbb R^{p_2}_-$ and $h:\mathbb R^n\times\mathbb R^m\rightarrow\mathbb R^{p_1}, g:\mathbb R^n\times\mathbb R^m\rightarrow\mathbb R^{p_2}$ are continuously differentiable. Since the restricted inf-compactness is weaker than the uniform compactness of the solution map near $\bar x$, the directional FOSCMS is weaker than MFCQ, and the directional solution set is smaller than the solution set,  our results for the case of parametric nonlinear program have improved the classical result of Gauvin and Dubeau \cite[Corollary 4.3]{GD}. 
 
\begin{thm} Let $u\in \mathbb{R}^n$ be a given direction. Assume that for the parametric nonlinear program, the restricted inf-compactness holds at $\bar x$ in direction $u$ with compact set $\Omega_u$. Suppose for each $y\in S(\bar x)\cap\Omega_u$ and  every $(d^h,d^g)\in{\rm Im}D_xP(\bar x, y;u)$, there exists $v\in \mathbb{R}^m$ satisfying
		\begin{align}
		&d^h+\nabla_y h(\bar x, y)v=0, \label{eqn37}\\
		&(d^g+\nabla_y g(\bar x, y)v)_j\leq0\ \mbox{for any}\ j\in\mathcal I_g(\bar x, y):=\{j=1,\ldots,p_2|g_j(\bar x, y)=0\}. \label{eqn38}
		\end{align}
	Suppose that FOSCMS of the system $h(\bar x,y)=0, g(\bar x,y)\leq0$, i.e,
	\begin{equation*}
		\nabla_y h(\bar x, y)^T\mu+\nabla_y g(\bar x, y)^T\gamma=0,\ 0\leq\gamma\perp g(\bar x, y), \gamma\perp  w\implies(\mu,\gamma)=(0,0)
	\end{equation*} holds
	for every direction $w\in D(\bar x, y;u)$ and  each $y\in {S(\bar x)\cap\Omega_u}$,
	where
	\begin{eqnarray*}
\lefteqn{	D( \bar x, y;u)} \\
& := & \left \{ \left  (d^h+\nabla_y h( \bar x, y)v, d^g+\nabla_y g( \bar x, y)v \right )\mid \begin{array}{l}
	0\not =(d^h,d^g,v) \in {\rm Im} D_x P(\bar  x, y;u)\times \mathbb{R}^m\\
(d^h,d^g)	\mbox{ satisfies (\ref{eqn37})-(\ref{eqn38})} 
	\end{array} \right \}.
	\end{eqnarray*}
	Then  RS holds for the problem $(P_{\bar x})$ 
	at each $y\in {S(\bar x)\cap\Omega_u}$ in direction $u$.
	Suppose  for  each $y\in {S(\bar x)\cap\Omega_u}$ and each $v$  satisfying (\ref{eqn37})-(\ref{eqn38}) with $(d^h,d^g)\in{\rm Im}D_xP(\bar x, y;\pm u)$
	\begin{eqnarray*}
	\left \{ \begin{array}{l} 	\nabla_y h(\bar x, y)^T\mu+\nabla_y g(\bar x, y)^T\gamma=0,\ 0\leq\gamma\perp g(\bar x, y),\\ 
	  \mu\perp  \nabla h(\bar x,y)(\pm u,v),
	 \gamma\perp  \nabla g(\bar x,y)(\pm u,v) \end{array} \right .\Longrightarrow  (\mu,\gamma) =(0,0).
	\end{eqnarray*}
	Then \begin{align*}
		{\min_{ y\in S(\bar x;u)\cap\Omega_u}}
		{\min_{\lambda\in  \Lambda(\bar x, y)}}\nabla_x L(\bar x,y, \lambda)u\leq V_-'(\bar x;u) \leq V_+'(\bar x;u)\leq{\min_{ y\in S(\bar x;u)\cap\Omega_u}}
		{\max_{\lambda\in  \Lambda(\bar x, y)}}\nabla_x L(\bar x,y, \lambda)u.
	\end{align*} 
	Furthermore if  $\Lambda(\bar x, y)=\{\bar \lambda(y)\}$ is a singleton for every 
	$y\in S(\bar x;u)\cap\Omega_u$, 
	then $V(x)$ is Hadamard directionally differentiable and 
	$$V'(\bar x; u)= \min_{ y\in S(\bar x;u)\cap\Omega_u}
	\nabla_x L(\bar x,y, \bar \lambda(y))u.$$
\end{thm}
\beginproof 
Since $C$ is a convex polyhedral set, by Remark \ref{Remark3.1}, condition (\ref{geom}) in Proposition \ref{SRS} is reduced to the feasibility of the system (\ref{assumpv}). Using the expression of the tangent cone the system (\ref{geom}) is (\ref{eqn37})-(\ref{eqn38}). The conclusion of RS in direction $u$ then follows from Proposition \ref{SRS}.
By Remark \ref{wRCQl}, since $C$ is a convex polyhedral set, RS implies the feasibility of system (\ref{feasibleso}). After translating the regularity condition (\ref{wRCQ4}) to our setting, we obtain the rest of results from Theorem \ref{nondir}.
\endproof
Similarly, when the solution map $S$ enjoys  the directional inner semicontinuity,  by Corollary \ref{isc2} 
one has the following  results with easier-to-verify assumptions.
\begin{cor}\label{Cor4.7} Let $u\in \mathbb{R}^n$ be a given direction. Suppose that for the parametric nonlinear program, the solution map $S(x)$ is inner semicontinuous at $(\bar x,\bar y)\in{\rm gph}S$ in direction $u$. Then $\bar y\in S(\bar x;u)$. Suppose 
$(\nabla_xh(\bar x,\bar y)u,\nabla_xg(\bar x,\bar y)u)\neq(0,0)$,  there exists
$ (\alpha,v)\in\mathbb R\times\mathbb R^m$ such that 
\begin{equation} \nabla h(\bar x,\bar y)(\alpha u,v)=0, \nabla g_j(\bar x,\bar y)(\alpha u,v)\leq0, j\in\mathcal I_g(\bar x,\bar y),  \label{eqn39}\end{equation}  
	and  for every $(0,0)\neq(\alpha,v)$ satisfying (\ref{eqn39}), one has
	\begin{eqnarray*}
		 \left \{ \begin{array}{l}
		\nabla_yh(\bar x,\bar y)^T\mu+\nabla_yg(\bar x,\bar y)^T\gamma=0,0\leq\gamma\perp g(\bar x,\bar y), \\
		 \mu\perp  \nabla h(\bar x,y)(\alpha u,v),  \gamma\perp \nabla g(\bar x,\bar y)(\alpha u,v)
		 \end{array}
		 \right .
		&&  \implies (\mu,\gamma)=(0,0).
	\end{eqnarray*} 
	Then \begin{align*}
		{\min_{\lambda\in  \Lambda(\bar x,\bar y)}}\nabla_x L(\bar x, \bar y, (\mu,\gamma))u\leq V_-'(\bar x;u) \leq V_+'(\bar x;u)\leq
		{\max_{\lambda\in  \Lambda(\bar x,\bar y)}}\nabla_x L(\bar x, \bar y, (\mu,\gamma))u.
	\end{align*} 
	Furthermore if  $\Lambda(\bar x,\bar y)=\{(\bar\mu,\bar\gamma)\}$ is a singleton, 
	then $V(x)$ is Hadamard directionally differentiable and 
	$$V'(\bar x; u)=
	\nabla_x L(\bar x,\bar y, (\bar\mu,\bar\gamma))u.$$
\end{cor}
}

By Theorem \ref{cor4.3} and Corollary \ref{cor-nodir}, since the set $C$ is convex polyhedral for the parametric nonlinear program, we also have the following corollary.
 
\begin{cor} \label{cor4.10}  Consider the affine case of the parametric nonlinear program where  $$(h(x,y),g(x,y))=a(x)+By+c$$ with $a:\mathbb R^n\rightarrow\mathbb R^p$ being continuously differentiable, $B\in\mathbb R^{p\times m}, c\in \mathbb{R}$,  and  {${\cal F}(x)$ is nonempty in a neighborhood of $\bar x$ in direction $u$}. Assume that the restricted inf-compactness holds at $\bar x$ in direction $u$ with compact set $\Omega_u$.
Then 
$$ {\min_{ y\in S(\bar x;u)\cap\Omega_u}}
 {\min_{\lambda\in  \Lambda(\bar x, y)}}\nabla_x L(\bar x,y, \lambda)u\leq V_-'(\bar x;u) \leq V_+'(\bar x;u)
	\leq{\min_{ y\in S(\bar x;u)\cap\Omega_u}}
	 {\max_{\lambda\in  \Lambda(\bar x, y)}}\nabla_x L(\bar x,y, \lambda)u.$$
	 Furthermore  if  $\Lambda(\bar x, y)=\{\bar \lambda(y)\}$ is a singleton for every 
$y\in S(\bar x;u)\cap\Omega_u$,
then $V(x)$ is Hadamard directionally differentiable  at $\bar x$ in direction $u$ and 
$$V'(\bar x; u)= {\min_{ y\in S(\bar x;u)\cap\Omega_u}}
\nabla_x L(\bar x, y, \bar \lambda(y))u.$$

Alternatively, assume that $S(x)$ is inner semicontinuous at $(\bar x,\bar y)\in {\rm gph} S$ in direction $u$. Then $\bar y\in S(\bar x;u)$ and 
$$ 
 \min_{\lambda\in  \Lambda(\bar x, \bar y)}\nabla_x L(\bar x,\bar y, \lambda)u\leq V_-'(\bar x;u) \leq V_+'(\bar x;u)
	\leq 
	 \max_{\lambda\in  \Lambda(\bar x, \bar y)}\nabla_x L(\bar x,\bar y, \lambda)u.$$
	 Furthermore  if  $\Lambda(\bar x, \bar y)=\{\bar \lambda \}$ is a singleton,
then $V(x)$ is Hadamard directionally differentiable  at $\bar x$ in direction $u$ and 
$V'(\bar x; u)= 
\nabla_x L(\bar x, \bar y, \bar \lambda)u.$
\end{cor}
 The following example shows the advantage of using the directional solution map. 
\begin{example}
	Consider the following parametric nonlinear program
	\begin{align*}
		{\min_{y}f(x,y):=xy}\quad
		s.t.\ -y-x^2-1\leq0, y-x^2-1\leq 0.
	\end{align*}
Its solution map is
\begin{eqnarray*}
	S(x)=\left\{
	\begin{array}{ll}
	-x^2-1,\ &x>0,\\
	 \mbox{$[-1,1]$},\ &x=0,\\
	x^2+1,\ &x<0.
	\end{array}
	\right.
\end{eqnarray*}Let $\bar x=0$. Then $S(\bar x)=[-1,1]$. 
The Lagrange function is
$$L(x,y,(\lambda_1,\lambda_2))=xy+ (-y-x^2-1)\lambda_1+(y-x^2-1)\lambda_2 .$$
Linear independence constraint qualification (LICQ) holds at each $y\in S(\bar x)$. The set of multipliers $\Lambda(\bar x,  y)=\{(0,0)\}$ is  singleton for any $y\in [-1,1]$.  Since LICQ holds at each $y\in S(\bar x)$, one can use the classical result (\ref{minmax}) to obtain the directional differentiability of $V$ at $\bar x$ in direction $u$:
 $$V'(\bar x; u)=\min_{y\in S(\bar x)} L_x (\bar x, y, (0,0))u=\min_{y\in [-1,1]} yu.$$ Note that  there is a minimization to perform in order to obtain the directional derivative. We now apply our result. 
 $S(\bar x;1)=\{-1\}\subsetneq S(\bar x)$ and $S(\bar x;-1)=\{1\}\subsetneq S(\bar x)$.
The solution map  $S(x)$ is not inner semicontinuous at $(0,1)$ or 
$(0,-1)$, but is obviously inner semicontinuous at $(0,-1)$ and $(0,1)$ in direction $1$ and $-1$, respectively. Note that ${\cal F}( x)=[-x^2-1,x^2+1]$ is nonempty.
 Hence by either Corollary \ref{Cor4.7} or Corollary \ref{cor4.10}, $V(x)$ is Hadamard directionally differentiable at $\bar x$ in  direction $u=\pm 1$. 
 Since $S(\bar x;1)=\{-1\}$ and $S(\bar x;-1)=\{1\}$,
 we have $$V'(\bar x;1)=\nabla_x L(\bar x,-1, (0,0))(1)=-1,\ V'(\bar x;-1)=\nabla_x L(\bar x,1, (0,0))(-1)=-1.$$
\end{example}


\end{document}